\documentclass[authoryear,review]{elsarticle}
\usepackage{array}
\usepackage{breqn}
\usepackage{subcaption}
\usepackage{stackengine}
\usepackage{verbatim} 
\newcolumntype{R}{>{\centering\arraybackslash}m{8cm}}
\usepackage{lineno,hyperref}
\modulolinenumbers[5]
\journal{Journal of Computers and Operations Research}
\usepackage{array}
\usepackage{subcaption}
\usepackage{nomencl}
\usepackage{amsmath}
\usepackage{algorithm}
\usepackage{algorithmic}
\makenomenclature
\usepackage{etoolbox}
\renewcommand\nomgroup[1]{%
  \item[\bfseries
  \ifstrequal{#1}{A}{Sets}{%
  \ifstrequal{#1}{B}{Parameters}{%
  \ifstrequal{#1}{C}{Decision Variables}{%
  \ifstrequal{#1}{D}{Functions}{%
  \ifstrequal{#1}{E}{Problems}{}}}}}%
  ]}

\begin{document}

\begin{frontmatter}

\title{Test Center Location Problem: A bi-objective Model and Algorithms}

\author[a,b]{Mansoor Davoodi\corref{cor1}}
\ead{m.davoodi-monfared@hzdr.de}

\author[a,b,c,d]{Justin M. Calabrese}
\ead{j.calabrese@hzdr.de}

\cortext[cor1]{Corresponding author.}

\address[a]{Center for Advanced Systems Understanding, Untermarkt 20, 02826 Görlitz}

\address[b]{Helmholtz-Zentrum Dresden-Rossendorf, Bautzner Landstraße 400, 01328 Dresden}

\address[c] {Department of Ecological Modelling, Helmholtz Centre for Environmental Research – UFZ, Leipzig, Germany}

\address[d]{Department of Biology, University of Maryland, College Park, MD, USA}

\begin{abstract}
The optimal placement of healthcare facilities, including the placement of diagnostic test centers, plays a pivotal role in ensuring efficient and equitable access to healthcare services.  However, the emergence of unique complexities in the context of a pandemic, exemplified by the COVID-19 crisis, has necessitated the development of customized solutions. This paper introduces a bi-objective integer linear programming model designed to achieve two key objectives: minimizing average travel time for individuals visiting testing centers and maximizing an equitable workload distribution among testing centers. To address this problem, we propose a customized local search algorithm based on the Voronoi diagram. Additionally, we employ an $\epsilon$-constraint approach, which leverages the Gurobi solver. We rigorously examine the effectiveness of the model and the algorithms through numerical experiments and demonstrate their capability to identify Pareto-optimal solutions. We show that while the Gurobi performs efficiently in small-size instances, our proposed algorithm outperforms it in large-size instances of the problem.

\end{abstract}

\begin{keyword}
Testing center \sep Facility location \sep \textit{k}-balance \sep \textit{k}-median \sep Bi-objective optimization \sep Heuristics
\end{keyword}

\end{frontmatter}

\section{Introduction}

The purposeful allocation of facilities, which includes the selection of examination centers, has undergone thorough scrutiny spanning various fields, including operations research, geography, and transportation planning. Particularly, facility location (FL) problems within the context of healthcare systems have garnered significant attention due to their practical implications in enhancing healthcare delivery. These problems involve determining strategic locations for healthcare facilities such as hospitals, clinics, and medical centers to serve a given population while considering various factors, including geographical distribution of population, patient demand, resource constraints, and cost considerations. 
One of the recent challenges in this area was finding optimal locations for testing centers during the COVID-19 pandemic. 

As observed by many people worldwide in recent times, crowded conditions not only prolong waiting times in testing queues for individuals but also lead to increased chances of viral transmission. Hence, when increasing the number of test centers is not possible (e.g., due to resource constraints), the strategic organization and placement of these facilities assume paramount importance. This problem entails a delicate trade-off between two essential objectives. Firstly, it is imperative to minimize the traveling distance between individuals seeking testing and their nearest testing center, thereby reducing their associated travel time costs. On the other hand, in order to mitigate the risk of infection and ensure efficient service delivery, it is equally crucial to achieve an equitable workload distribution across these centers. Striking this balance should significantly contribute to the effective management of these facilities, fostering a fair distribution of responsibilities among them.
In light of two fundamental real-world considerations, namely, the tendency of individuals to select the closest available center and the constraints imposed by a limited number of such testing centers, the overarching objective became clear. The goal is to strategically deploy a limited (say $k$ centers) set of test centers with the dual purpose of minimizing the distance between individuals and their closest test center, while concurrently minimizing the differences in workload among the test centers.


In this paper, we consider clusters of individuals as weighted demand points, analogous to the population residing in an apartment complex, in conjunction with a predefined set of potential locations for establishing test centers. Consequently, the problem at hand entails the selection of \textit{k} potential center locations, wherein two key objectives are pursued: (\textit{i}) the attainment of maximum balance among the workloads of the centers, specifically minimizing the disparity between the highest and lowest workloads, and (\textit{ii}) the minimization of the average traveling time for the demand points. We name this particular FL problem the \textit{test center location problem} and denote it by TCLP.

It is pertinent to note that solving either of the objectives in the TCLP is an NP-hard problem. The first objective, often referred to as the 'k-balanced' objective, has been studied recently, while the second objective aligns with the well-established 'k-median' problem, which has received substantial research focus over time. In this study, we undertake a comprehensive approach to address this bi-objective optimization challenge. Initially, we formulate the problem as an integer linear program, providing a solid foundation for subsequent analysis.
We then proceed to propose two distinct approaches for obtaining Pareto-optimal solutions. The first approach involves leveraging the $\epsilon$-constraint approach in conjunction with the commercial solver Gurobi. This approach demonstrates efficacy, particularly for smaller problem instances; however, it exhibits notable computational demands for larger-scale scenarios. Consequently, as a second approach, we introduce a custom-designed bi-objective hill-climbing strategy that leverages geometric information such as the Voronoi diagram. Our implementation and comparative evaluation of these two approaches encompass a variety of problem instances, considering criteria such as runtime efficiency and the ability to identify Pareto-optimal solutions.
The simulation results highlight the superior performance of the proposed heuristic approach, underscoring its potential as a valuable tool for addressing the intricate challenges inherent in this bi-objective FL problem.

This paper consists of six sections. Section \ref{relatedwork} reviews prior research in FL, with a specific focus on healthcare facilities. Section \ref{modeling} formulates the TCLP and presents the integer linear program. Section \ref{approach} proposes the $\epsilon$-constraint method using the Gurobi solver as well as a bi-objective hill-climbing approach for solving the TCLP. Section \ref{simulation} discusses simulation results and provides a comparative analysis. Finally, Section \ref{conclusion} concludes the paper and outlines future research directions.



\section{Related work}
\label{relatedwork}
The FL problem requires the determination of appropriate locations (centers or hubs) of a set of facilities among a set of demand points (customers or clients) \cite{daskin1995network}. This problem has numerous real-world applications and has been widely studied in the literature of operations research, industrial engineering, applied mathematics, and computer science \cite{daskin1995network, farahani2010multiple, kochetov2005computationally}. There are several parameters, constraints, and objectives in the FL problem, and consequently, many variations of it have been studied \cite{daskin1995network, farahani2010multiple}. For example, the demand set may be discrete or continuous, weighted or unweighted, static or dynamic, certain or uncertain. The potential facility set can be discrete or continuous, and capacitated or incapacitated. Furthermore, several definitions for the objective function have been considered \cite{daskin1995network, farahani2010multiple, megiddo1984complexity}.

The objective function in FL problems, which is usually determined with regard to the type of application, is very important in the complexity class of the problems \cite{daskin1995network, farahani2010multiple, kochetov2005computationally}. For example, \textit{k-median} and \textit{k-center} are two well-known types of FL problems for public FL and emergency FL with the objectives \textit{min-sum} and \textit{min-max}, respectively. The NP-completeness of both of the problems (and some variations of them) has been proved \cite{megiddo1984complexity}, and many approximations and heuristic approaches have been proposed for solving them (e.g. see \cite{vazirani2013approx, davoodi2011solving, drezner1984p, mahdian2006approximation}).

In both \textit{k-median} and \textit{k-center} problems, the goal is optimizing the process for the client side, e.g., minimizing the average and maximum distance of each client from its closest center, which is useful for both public and emergency facilities. Both of these objectives belong to the client side, that is, objectives to emphasize the service quality that the clients receive. However, there are objectives such as the recently proposed \textit{k-balanced} objective that enhance the quality or eligibility in the center side \cite{davoodi2023bi}.
The \textit{k}-balanced objective focuses on the fair distribution of accessibility among the clients \cite{bortnikov2012load}. For example, consider the problem of placing some congruent antennas in a wireless network \cite{kleinberg1999fairness}. For some technical reasons, and to have a good connection quality, usually each client is assigned to its closest antenna(s). Thus, to manage the traffic in the network, it is necessary for the antennas to have almost the same network load. As another example, assume locating $k$ voting stations under the assumption that each person goes to their closest voting station. So to balance the crowding in the stations, the stations' workloads need to be balanced. These considerations may also apply in placing banks, stores, educational, cultural, and sports centers, and are very important in Territory Design \cite{kalcsics2005towards}. Note that in the \textit{k}-balanced problem, each client is served by the closest center; consequently, it is not an assignment problem \cite{bortnikov2012load}.

Similar to the FL problems, there are many parameters and constraints in the \textit{k}-balance problem, and different variations of it can be presented. In addition to the discrete or continuous potential facility centers and different metrics, the definition of the term ``maximum balance'' is not unique and can be determined by the type of application. Marín \cite{marin2011discrete} originally proposed the \textit{k}-balance problem in 2011. He studied the discrete version of the problem and constructed integer programming formulations of a variation of the problem and proposed a branch-and-cut algorithm for solving them \cite{marin2011discrete}. Finally, he evaluated the algorithm by some simulations that used computational time as the efficiency factor \cite{marin2011discrete}. He noted that the number of valid inequalities in the formulations of the problem is exponential. Filipović et al. \cite{filipovic2012modification} proposed a combined heuristic method consisting of a genetic algorithm with an interchange heuristic for the balanced allocation problem \cite{filipovic2012modification}. This combined method was a variable neighborhood search heuristic that utilizes a technique called shaking neighborhood in order to avoid becoming stuck in local optima, which has subsequently been improved by Kratica et al. \cite{kratica2012variable}. 
Davoodi \cite{davoodi2019k} originally discussed the complexity of the \textit{k}-balance problem with two different objectives: (\textit{i}) minimizing the maximum number of allocated clients to any center, and (\textit{ii}) minimizing the difference between the maximum and the minimum number of clients allocated to the center. He showed NP-hardness of the \textit{k}-balance problem for both objectives in the plane under both Manhattan and Euclidean metrics.

FL in the context of healthcare is a multifaceted challenge that involves optimizing accessibility, resource allocation, and patient outcomes. The related work in this field encompasses various modeling techniques, GIS applications, and patient-centric approaches to address the complex task of facility placement \cite{daskin2004location, ahmadi2017survey}. 
Flores et al. \cite{flores2021optimizing} focused on healthcare FL in low and middle-income countries, particularly the Philippines. They introduced a novel cooperative covering maximal model to optimize primary care facility placement using open-source data, considering equity and efficiency parameters. The approach holds promise for evidence-based healthcare facility decisions in resource-limited settings and can be adapted to other sectors. 

Liu et al. \cite{liu2023location} aimed to explore the principles and factors impacting the choice of locations for emergency medical facilities during public health crises. They delved into the process of identifying optimal facilities and introduced a logistic regression model to establish a site selection framework tailored for emergency medical facilities in megacities during public health emergencies. 
Karmel et al. \cite{shehadeh2021equity} addressed equity in stochastic healthcare FL models, examining how uncertainty affects disparities. They focused on modeling uncertainty, equity, and FL, encompassing aspects and outcomes like tractability, fairness, and access metrics.
Wang et al. \cite{wang2018healthcare} studied the FL problem in China's evolving healthcare landscape, particularly, location-allocation challenges in growing cities. They introduced a hierarchical model balancing social, economic, and environmental factors, using a bi-level multi-objective particle swarm optimization algorithm for complex decisions.
Fathollahi et al. \cite{fathollahi2021multi} addressed the global challenge of an aging population, whereby healthcare decision-makers face the complexities of optimizing home healthcare for the elderly and ensuring its sustainability. They introduced a robust multi-objective optimization model for home healthcare, considering factors like caregiver scheduling, care continuity, patient availability, service times, and quality standards. Finally, they presented a metaheuristic to tackle the problem.

Tang et al. \cite{tang2022bi} studied a multi-period vaccination planning problem, optimizing vaccination recipients' travel distance and operational costs. The problem involves deciding when to open vaccination sites, how many stations to launch, recipient assignments, and site replenishment. Initially framed as a bi-objective mixed integer linear program, they introduced a weighted-sum, $\epsilon$-constraint and used genetic algorithms to solve the problem.
Alhothali et al. \cite{alhothali2022location} discussed the COVID-19 vaccination center location problem with the objectives of enhancing accessibility and minimizing costs. They employed maximal coverage models with a focus on minimizing transportation time and travel distance. Maliki et al. \cite{maliki2022multi} studied multi-period FL decisions in scenarios emphasizing pandemics with volatile demand, and including opening, relocating, closing, and utilizing mobile facilities. They employed NSGA-II to balance economic costs and CO2 emissions. 
Lai et al. \cite{lai2021multi} presented a vaccination station location model, incorporating multi-period planning for medical professionals, vaccine procurement, and inventory decisions amidst demand uncertainties. Formulated as a complex two-stage stochastic problem, they utilized a Benders decomposition-based heuristic for effective resolution.


\section{Test center location problem}
\label{modeling}
The test center location problem (TCLP) seeks to determine the optimal arrangement of a set of $k$ test centers in a manner that simultaneously minimizes the travel cost for individuals (who are typically tested by their closest center, reflecting real-world conditions) and maximizes the equitable distribution of workload among these centers. Given that we are examining identical test centers, we define a center's workload as the count of individuals it serves. Additionally, we take into account the average travel time between an individual and their closest center. We define the objective of workload balance as the minimization of the disparity between the most heavily populated center and the least heavily populated one. Since individuals are tested in the nearest center, there exists a trade-off between travel time and the balance of workload among the centers. Solving such a problem provides a set of \textit{Pareto-optimal} solutions, i.e. those that cannot be enhanced in one objective without compromising the other objective.
Within this section, we formulate the TCLP formally, and subsequently, we articulate an integer linear programming model tailored to address the problem. 

\subsection{Test center location problem formulation}
To establish a comprehensive model test center location problem, we introduce the concept of assigning weights to each demand point, with each weight corresponding to the population count residing at that particular demand point. This weighted approach significantly enhances the problem's applicability to real-world scenarios. For instance, all individuals residing in an apartment complex can be represented as a single demand point, with its weight equal to the number of residents within it. In larger-scale instances of the problem, such as those involving extensive urban areas, it becomes feasible to preprocess the data by clustering residents who are in proximity. The center of each cluster is then assigned a weight equivalent to the size of that cluster. This preprocessing step leads to a substantial reduction in the problem's dimensionality, ultimately facilitating the proposal of an efficient solution.
We now define the notation and problem formulation precisely.

\printnomenclature[1cm]


\nomenclature[B, 01]{$n$}{Number of demand points}
\nomenclature[B, 02]{$P=\{(p_1,w_1),(p_2,w_2),...,(p_n,w_n)\}$}{ Set of weighted demand points}
\nomenclature[B, 03]{$p_i= (x_i,y_i)$}{$i$-th demand point with coordination of $(x_i,y_i)$ in the plane}
\nomenclature[B, 04]{$w_i$}{Weight of $i$-th demand point, that is, number of individuals located in location $p_i$}
\nomenclature[B, 05]{$m$}{Number of potential test center locations}
\nomenclature[B, 06]{$Q=\{q_1,q_2,...,q_m\}$}{ Set of potential test center locations}
\nomenclature[B, 07]{$q_j=(x_j,y_j)$}{$j$-th potential test center which is located in coordinate $(x_j,y_j)$}
\nomenclature[B, 08]{$d_{ij}$}{Traveling distance (or any type of cost in general) between $p_i$ and $q_j$}
\nomenclature[B, 09]{$k$}{Number of test centers which must be chosen (or say to be \textit{opened})}
\nomenclature[B, 10]{$C=\{c_1,c_2,...,c_k)\}$}{Set of opened centers. $C \subseteq Q$, $C$ is a (feasible) solution}
\nomenclature[B, 11]{$\delta(p_i)$}{ Closest opened center to $p_i, \forall i \in \{1,2,...,n\},~~\delta(p_i) \in C$}
\nomenclature[B, 12]{$\Delta(c_j)$}{ All demand points whose closest opened center is $c_j$, ($\Delta(c_j)=\{p_i \in P : c_j = \delta(p_i)\}$)}
\nomenclature[B, 13]{$u$}{Maximum number of (weighted) demand points allocated to any opened center, $u = \max \limits_{c \in C} \sum \limits_{p_i \in \Delta(c)} w_i$}
\nomenclature[B, 14]{$l$}{Minimum number of (weighted) demand points allocated to any opened center $l = \min \limits_{c \in C} \sum \limits_{p_i \in \Delta(c)} w_i$}

Given a weighted set of demand points, denoted as $P = \{(p_1, w_1), (p_2, w_2), \ldots, (p_n, w_n)\}$, a set of potential facility centers $Q = \{q_1, q_2, \ldots, q_m\}$, the travel distance $d_{ij}$ for any pair $(p_i,q_j)$, and an integer $k$, the goal is to select (or open) $k$ centers from the available $m$ potential locations such that achieves the following objectives:

The first objective function, workload balance, is defined as follows
\begin{equation}
Workload~~balance:~~ Minimize~~ F_1(C) =  u - l, 
\label{F1_Obj}
\end{equation}

where $u$ and $l$ are the maximum and minimum number of individuals (demand points) allocated to any opened center, respectively.
The second objective function is \textit{k-median} objective, that is
\begin{equation}
Weighted~min~sum:~~Minimize~~ F_2(C) =  \frac{1}{\Sigma_{i=1}^{n} w_i}  \Sigma_{i=1}^{n} w_i d_{i\delta(p_i)}~,
\label{F2_Obj}
\end{equation}

where $d_{i\delta(p_i)}$ denotes the traveling distance between demand point $p_i$ and its closest opened center, $\delta(p_i)$. We assume $\delta(p_i)$ is unique. Therefore, this objective aims to minimize the average (weighted) travel distance for all individuals. It's important to note that there are no stringent constraints other than the requirement to precisely open $k$ centers from the initial set of $m$ potential centers, a choice usually influenced by financial or expertise limitations, such as constraints on available nurses or doctors.

\subsection{Integer linear programming model for the test center problem}
In the following sections, we provide a bi-objective Integer Linear Programming (ILP) model. The model is based on the formulation presented by Marín \cite{marin2011discrete}, which we extend for the weighted demand points and the two contrasting objectives. To this end, we define the following binary variables:

\begin{equation}
x_{ij} = 
\begin{cases}
1 & \text{, if demand point $p_i$ is served by center $c_j$ (or if $p_i \in \Delta(c_j)$)} \\
0 & \text{, otherwise}
\end{cases}
\end{equation}

\begin{equation}
 y_{j} = 
\begin{cases}
1 & \text{, if center $c_j$ is selected to be opened (or if $c_j \in C$)} \\
0 & \text{, otherwise}
\end{cases}
\end{equation}

By having these $m(n+1)$ binary variables, the ILP for the test center location problem can be formulated as below:

\begin{equation}
\begin{aligned}
&Minimize~~ F_1 =  u - l\\
&Minimize~~ F_2 =  \frac{1}{\Sigma_{i=1}^{n} w_i}  \Sigma_{i=1}^{n} w_i (\Sigma_{j=1}^{m} x_{ij} d_{ij})\\
&~~~~Subject~to:\\
&~~~~~~~~~~~~\sum_{j=1}^{m} y_j = k\\
&~~~~~~~~~~~~~x_{ij} \leq y_j,~~\forall i \in \{1,2,...,n\},~~\forall j \in \{1,2,...,m\}\\
&~~~~~~~~~~~~\sum_{j=1}^{m} x_{ij} = 1~~\forall i \in \{1,2,...,n\}\\
&~~~~~~~~~~~~u \geq \sum_{i=1}^{n} w_i x_{ij}, ~~\forall j \in \{1,2,...,m\}\\
&~~~~~~~~~~~~l \leq \sum_{i=1}^{n} w_i x_{ij}+ (1-y_j) \sum_{j=1}^{m} w_{i},~~\forall j \in \{1,2,...,m\}\\
&~~~~~~~~~~~~\sum_{j'=1}^{m} d_{ij'} x_{ij'}+ (M-d_{ij})y_j \leq M , ~~\forall i \in \{1,2,...,n\},~~\forall j \in \{1,2,...,m\}\\
&~~~~~~~~~~~~x_{ij},y_j\in \{0,1\},~~\forall i \in \{1,2,...,n\},~~\forall j \in \{1,2,...,m\}\\
\label{Model}
\end{aligned}
\end{equation}

The formulation is similar to the $k-median$ problem formulation, except the last constraint of the model, $\sum_{j'=1}^{m} d_{ij'} x_{ij'}+ (M-d_{ij})y_j \leq M$, where $M$ is a sufficiently large number (e.g., $M = \max \limits_{1\leq i \leq n,~ 1\leq j \leq m} d_{ij}$). This constraint guarantees each demand point is allocated to its closest center. We call this bi-objective ILP problem the \textit{test center location problem}, and denote it by TCLP.

To address the TCLP, an ideal algorithm should aim to yield a set of Pareto optimal solutions that exhibit diversity across the objective space. While the number of Pareto optimal solutions in the TCLP is finite, it can potentially be exponential in the worst case. To facilitate effective decision-making, the focus is to find a limited number of Pareto optimal solutions that cover all Pareto regions. This concept is commonly referred to as providing a "handful" of diverse Pareto optimal solutions \cite{deb2011multi, coello2007evolutionary}. Typically, this would encompass approximately 10 solutions, including not only extreme solutions for objectives $F_1$ and $F_2$ but also covering a broad spectrum of the objective space.
Then, the decision-maker can choose one of the provided Pareto-optimal solutions based on high-level information or any preferences that have not been integrated into the model. It is notable, that there are studies that suggest picking one Pareto-optimal solution like \textit{knee point} or other Nash solutions \cite{branke2004finding, gaudrie2018budgeted}. In the next section, we propose two approaches to find the Pareto optimal solutions of the TCLP.


\section{Solution approach for test center location problem}
\label{approach}
The test center location problem is classified as an NP-hard problem, signifying an absence of polynomial-time algorithms and rendering the task of identifying even a single Pareto optimal solution computationally infeasible. Consequently, the exploration of approximation and heuristic methods becomes invaluable. In this section, we first introduce an $\epsilon$-constraint approach capable of yielding a single Pareto optimal solution per execution. Furthermore, we propose a finely-tailored, efficient bi-objective hill-climbing approach designed to discover a set of \textit{non-dominated} solutions. Given the local-search nature of this approach, it is important to note that these non-dominated solutions may or may not represent the real Pareto optimal solutions. However, through extensive simulations and comparisons with the Pareto optimal solutions obtained via the $\epsilon$-constraint approach, we affirm that the majority of the non-dominated solutions either belong to the Pareto-optimal set or exhibit remarkable proximity to the Pareto-optimal fronts.

\subsection{An $\epsilon$-constraint method for the TCLP}
\label{epsiloncons}
The $\epsilon$-constraint method is a popular, simple and flexible method for multi-objective optimization, but it typically has limited ability to provide in-depth insights into Pareto-optimal solutions \cite{deb2011multi, chankong2008multiobjective}. In fact, the $\epsilon$-constraint method requires the designation of one objective as primary and the others as constraints. This categorization of primary and secondary of course can be subjective and may lead to biased results. So, the decision maker needs to have extra knowledge and perform additional analyses. After setting a proper upper threshold for the objective values that are designated as the constraints, this approach may efficiently work for well-distributed convex Pareto-optimal solutions but can be challenged by problems with densely concentrated search spaces.

In the implementation of the $\epsilon$-constraint method, we select $F_2$ and represent it as a constraint within the TCLP model, Eq.(\ref{Model}). To facilitate decision-maker comprehension and practicality, we establish both a lower and an upper bound for $\epsilon$ values. This ensures that $F_2$ objectives remain within this predefined range. Specifically, for this purpose, we set $k=m$ and compute $F_2$ values represented as $\frac{1}{\Sigma_{i=1}^{n} w_i}  \Sigma_{i=1}^{n} (w_i \min \limits_{1\leq j \leq m} d_{ij})$. This corresponds to solutions where each demand point is allocated to its closest potential facility center. Conversely, by configuring $k=1$, where all the demand points are allocated to the same center and $F_1$ is no longer important, the upper bound for the $F_2$ value can be established within polynomial time. Finally, the decision maker's preferences for the number of desired Pareto optimal solutions sets the number of $\epsilon$ values, which are then uniformly selected from this defined range and set in the following constraint.

\begin{equation}
\frac{1}{\Sigma_{i=1}^{n} w_i}  \Sigma_{i=1}^{n} w_i (\Sigma_{j=1}^{m} x_{ij} d_{ij}) \leq \epsilon,
\label{F2_cons}
\end{equation}

Upon relocating the previously mentioned equation to the constraint section, we transform the problem into a single-objective optimization model, i.e., as an ILP. This model can be efficiently solved using widely available commercial solvers like Gurobi ~\cite{gurobi}, yielding a single Pareto-optimal solution in each run. By introducing variations in the $\epsilon$ values and iteratively executing the process, we can systematically generate a diverse set of Pareto-optimal solutions.

\subsection{A local search approach for the TCLP}

One of the key factors in the success of heuristic and local search algorithms is the way that they generate a new solution using the obtained solutions in each iteration. The other factor is striking a balance between exploration and exploitation power. Considering these two factors, in this section, we propose a customized local search algorithm for the TCLP. We call this algorithm the \textit{Test Center Location Algorithm} and denote it by {TCLA}. This algorithm is a population-based heuristic algorithm that regenerates solutions by leveraging Voronoi neighbors. This regeneration method is called the \textit{Voronoi exchange operator} and ensures the gradual reproduction of new generations through an exchange operator, avoiding abrupt changes akin to the "mutation" process in genetic algorithms. Instead, it explores the search space in multiple directions facilitated by Voronoi analysis. In terms of exploitation, the population is updated using a "non-domination" comparison criterion. This entails retaining solutions that are non-dominated concerning the current population. In a bi-objective problem, a solution $C$ is said to \textit{dominate} a solution $C'$ if it is better in at least one objective and not worse in the other objective.
In cases where the number of such non-dominated solutions exceeds the population's capacity, a crowding operator such as a basic clustering technique is employed to select the most diverse non-dominated subset. The following sections will delve into the specifics of this process.

In a preprocess, we first compute Voronoi neighbors of each potential facility center, set $Q$. This can be performed in $O(m log m)$ time \cite{de1997computational}. Likewise, the nearest center to a given demand point can be determined in $O(\log m) time$. Let $Vor(q)$ denote Voronoi neighbors of each center $q\in Q$ in the Voronoi diagram. 
For a solution $C=\{c_1,c_2,…,c_k\} \subset Q$ for the TCLP, a random solution $C'$ can be generated by the following Voronoi exchange operator. We choose a center $c \in C$ and replace it with a random center $c' \in Vor(c)$. This Voronoi exchange operator can be applied for all centers, generating $k$ new random solutions. Since the evaluation process, which involves computing $F_1(C)$ and $F_2(C)$, is computationally expensive, we identify and remove repeated solutions before computing the objective values.

Now, let's elucidate the functioning of the entire local search algorithm. We commence with the assumption of a population denoted as $Pop_t$ (i.e., for t=0 in the beginning), possessing a size of $N$, and initialize it with random solutions like $C={c_1,c_2,…,c_k}$. Subsequently, we conduct an evaluation of these solutions, calculating their respective objective values, $F_1$ and $F_2$. This evaluation process demands $O(k \log k + n \log k)$ time for an individual solution $C$ by leveraging the Voronoi diagram of $C$ and identifying the closest center for each demand point.
In the next step, we employ the previously described Voronoi exchange operator to generate $k$ random solutions for every solution within the population, totaling $kN$ solutions in entirety. Following the removal of duplicated solutions and the computation of objective values for these generated solutions, we execute a non-dominated sorting, which identifies all non-dominated solutions within $O(kN \log (kN))$ time \cite{jensen2003reducing}.

In the final phase, we select non-dominated solutions from the union of $Pop_t$ and the newly generated solutions and construct $Pop_{t+1}$ with $N$ solutions. Two cases may happen, if the number of non-dominated solutions is less than $N$, we fill $Pop_{t+1}$ with the second level of non-dominated solutions. That is the non-dominated solutions after removing the first level. We repeat this process to fill $Pop_{t+1}$ with $N$ solutions. The second case happens if the number of non-dominated solutions exceeds $N$. In this case, we employ a \textit{crowding operator} to select a diverse ensemble of non-dominated solutions. Various approaches exist for reaching diversity among the solutions \cite{deb2011multi, coello2007evolutionary}. As an example, we first normalize the objective values and initiate by selecting two extreme solutions, those with the minimum $F_1$ and minimum $F_2$ values, incorporating them into $Pop_{t+1}$. Following this, we proceed to determine the largest axis-aligned bounding box that encompasses each solution while ensuring that no other solution is contained within it. We select the $N-2$ solutions that have the largest bounding boxes and incorporate them into $Pop_{t+1}$. This approach can be easily implemented by sorting the solutions based on their objectives. Consequently, it requires a time complexity of $O(Nk \log (kN))$.

Therefore, TCLA initiates its process with an initial random population and then proceeds to generate a new population through the utilization of the Voronoi exchange operator. From these populations, it selects the non-dominated solutions to be carried forward into the subsequent generation. These steps are reiterated for a specified number of iterations to accomplish its optimization objective. The pseudocode for TCLA is presented in \ref{Algo}. The time complexity of this algorithm for one iteration is $O(Nkn \log k)$ for evaluating the solutions using their corresponding Voronoi diagram, plus $O(Nk \log(Nk))$ if the crowding operator is needed.

It is worth noting that, we utilize the Voronoi diagram for a dual purpose: to ascertain the neighbors of a given solution and to expediently calculate the objective values associated with a solution. The number of Voronoi neighbors pertaining to a solution may exhibit variability, ranging from 2 to $(k-1)$. Nevertheless, the average number of Voronoi neighbors is constant. Additionally, the overall number of neighbors remains linear ($\leq 3k$). The assessment of a solution can be achieved through a brute-force algorithm in $O(nk)$ time; however, by employing the Voronoi diagram and performing the nearest point query, this process can be improved to $O(n \log k)$ time complexity \cite{de1997computational}.

TCLA, like all population-based heuristics, requires two predefined parameters: the size of the population ($N$) and the number of iterations. Remarkably, TCLA stands out by not requiring any additional parameters. In contrast, many heuristic algorithms necessitate a multitude of parameters, including crossover rate, mutation probability, and learning weights, among others \cite{coello2007evolutionary}. We firmly believe that in FL problems, particularly in the case of large instances, the Voronoi diagram plays a crucial role in efficiently achieving a balance between exploration and exploitation concepts within the search space. The Voronoi partition of the space serves as a valuable tool for distributing the combinatorial complexity of the problem into localized complexities.

\begin{algorithm} 
\caption{Test Center Location Algorithm (TCLA)}\label{Algo}
\textbf{Input:} Sets $P$ and $Q$, distance function (or matrix $d_{ij}$) and the integer number $k$\\
\textbf{Output:} Set of non-dominated solutions for the TCLP\\
\begin{algorithmic}[1]
\STATE Set the size of population to $N$, and number of generations to $T$
\STATE Initialize population $Pop_0$ with $N$ random solutions like $C=\{c_1,c_2,...,c_k\}$
\STATE For any solutions $C \in Pop_0$, compute Voronoi diagram of $C$, denoted by $VD(C)$, and then evaluate their objective values, $F_1(C)$ and $F_2(C)$
\FOR{$t=0$ \textbf{to} $T-1$}
\STATE For any solutions $C \in Pop_t$, apply $VD(C)$ and the Voronoi exchange operator and reproduce $k$ neighbor solutions. Put the new generated $Nk$ solutions in temp a population $TPop$
\STATE Remove the duplicated solutions in $TPop$
\STATE For any solution $C' \in TPop$, compute $VD(C')$, $F_1(C')$ and $F_2(C')$.
\STATE Add solutions in $Pop_t$ to $TPop$
\STATE Create an empty population $Pop_{t+1}$
\STATE Find all non-dominated solutions in $TPop$ and pop them into $Pop_{t+1}$
\IF {size of $Pop_{t+1} > N$}
\STATE Apply the crowding operator and choose $N$ most diverse non-dominated solutions.
\ELSE 
\WHILE {size of $Pop_{t+1} < N$}
\STATE Pop the non-dominated solutions from $TPop$ and add them to  $Pop_{t+1}$ if there exist some free slots, otherwise, put a random number of them to fill $Pop_{t+1}$ with $N$ solutions.
\ENDWHILE
\ENDIF
\ENDFOR
\STATE Return $Pop_T$

\end{algorithmic}
\end{algorithm}


\section{Simulation results}
\label{simulation}
This section is structured into two segments, presenting the outcomes of our proposed model and algorithms for identifying Pareto optimal solutions in the context of the TCLP. In the initial part, we employ the suggested TCLA on various problem instances with varying configurations. We present the outcomes both in the variable space and the objective space. In the subsequent part, we conduct a comparative analysis between TCLA and the $\epsilon$-constraint method, solved using the Gurobi solver. This comparison is made with regard to execution time and their respective capacities for identifying Pareto optimal solutions.

\textbf{Results on TCLA}

\begin{figure}[h]
  \centering
  \includegraphics[width=1.0\textwidth]{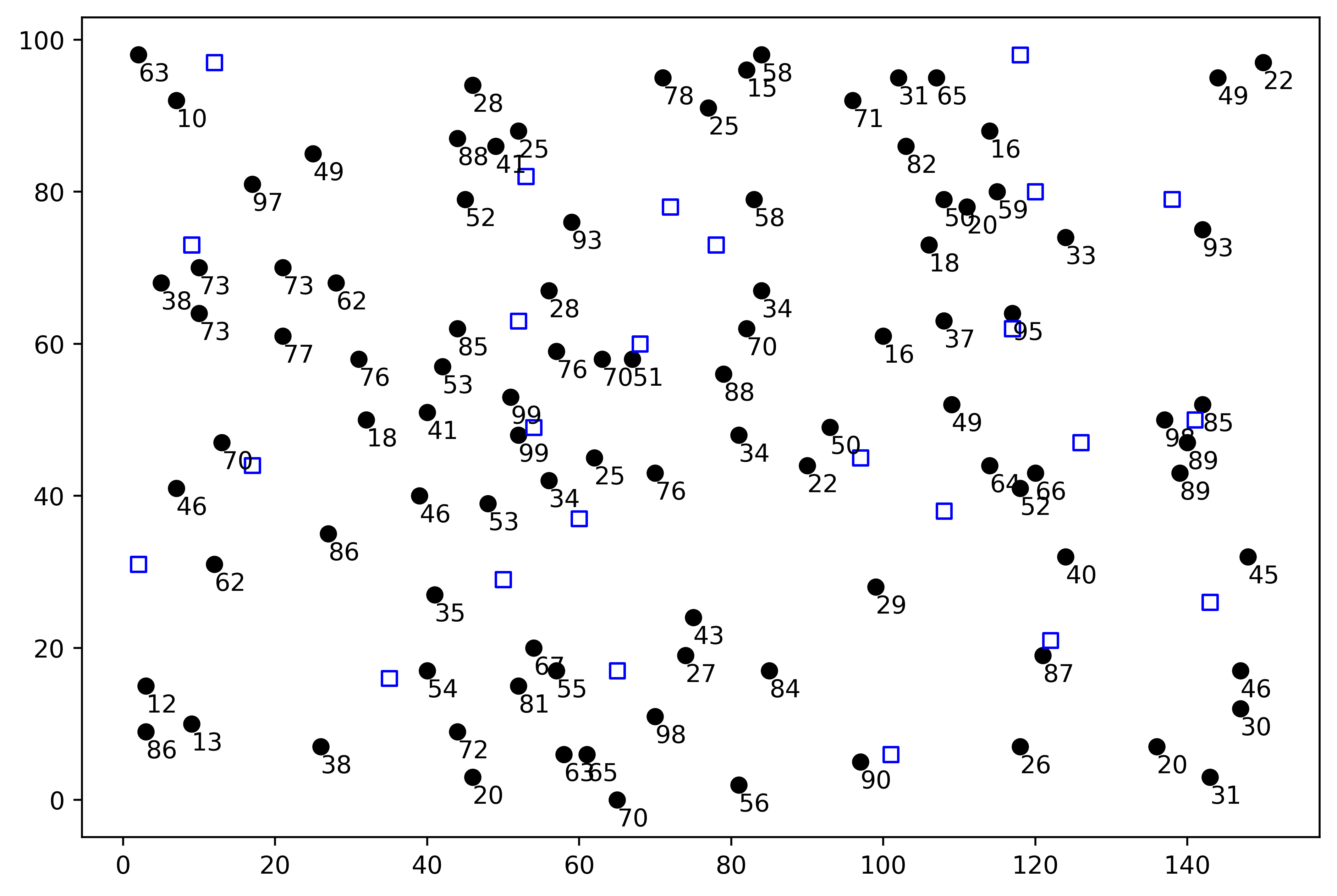}
  \caption{A random instance of the TCLP with $n=100$ weighted demand points (black circles) and $m=25$ potential test center locations (blue squares).}
  \label{randenv}
\end{figure}

We run TCLA on the model presented in Eq. (\ref{Model}) to find Pareto optimal solutions. The code of the algorithm is implemented in the programming language Python 3.7 and runs on a standard PC (\textit{Intel}(\textit{R}) \textit{Core}(\textit{TM}) \textit{i}7 and 32G \textit{RAM}). To this end, we consider a rectangular environment with a size of 1500x1000 and generate instances with random locations for the demand points, $P$, and potential centers, $Q$. Also, we assign random weights for the demand points in $[10,100]$. Figure \ref{randenv} shows the random instance with $n=100$ weighted demand points and $m=25$ locations at which test centers will be opened.

We denote each random instance of the TCLP with a triplet $(n,m,k)$, where $k$ is the number of opened centers. We run TCLA for $(100,25,k)$, where $k \in \{5,8,12,15\}$. The combinatorial complexity of the search space of the TCLP is related to $m$ and $k$ such that the number of possible solutions is $\binom{m}{k}$. This implies that the worst case happens for $k=\frac{m}{2}$. On the other hand, the complexity of TCLA, like all the population-based heuristics, is directly related to the size of the population, $N$, and the number of generations, $T$.
We set the size of the population in TCLA to $N = 2cm$ and the number of generations to $T = cN$, where $c=min\{k,m-k\}$. We choose these values because they achieve an optimal balance between processing time and the quality of the obtained non-dominated solutions. This choice is informed by our analysis of the simulation results.

\begin{figure}[h]
    \centering
    \begin{subfigure}{0.49\textwidth}
        \centering
        \includegraphics[width=\linewidth]{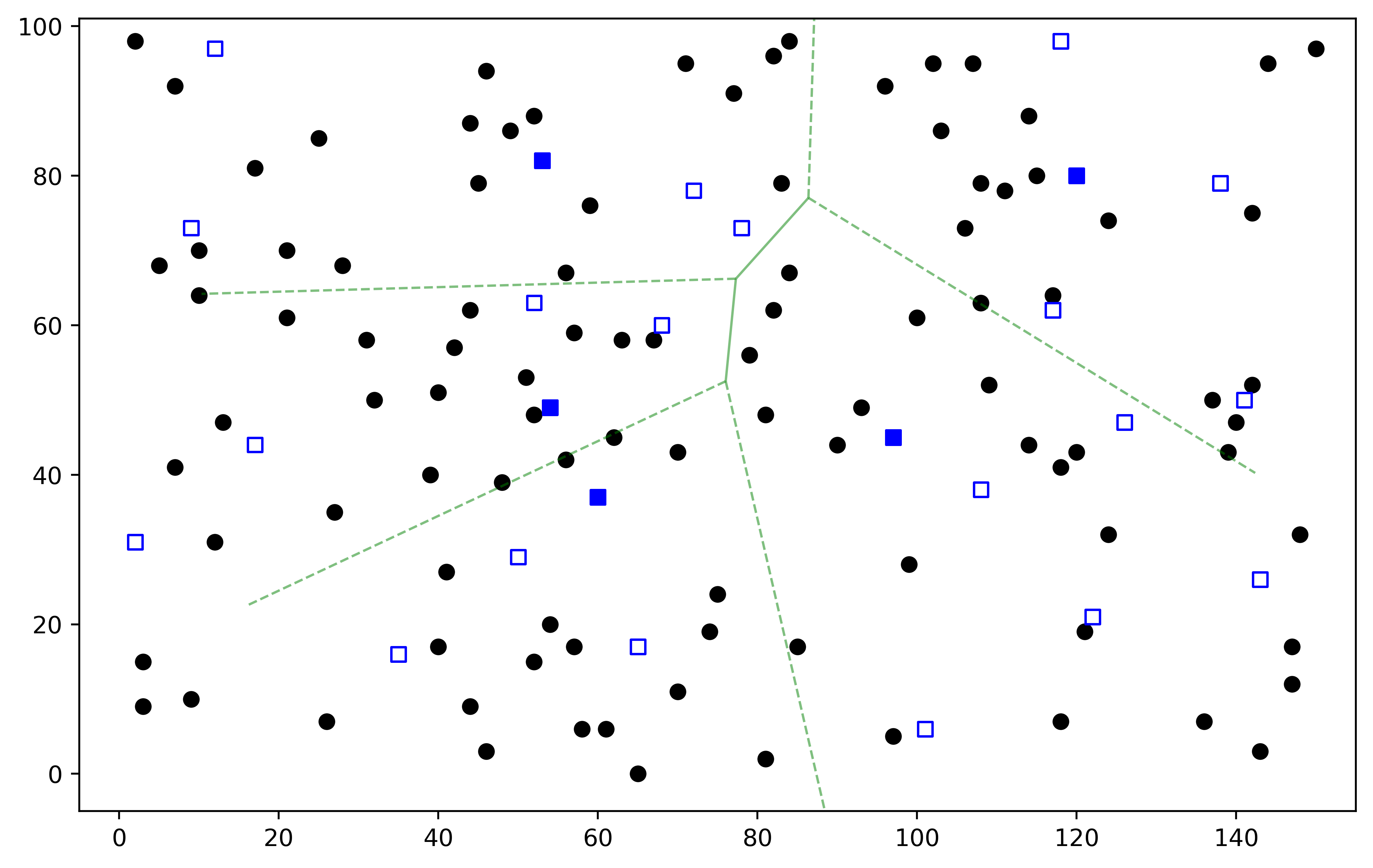}
        \caption{The obtained solution with minimum $F_1$}
    \end{subfigure}
    \hfill
    \begin{subfigure}{0.49\textwidth}
        \centering
        \includegraphics[width=\linewidth]{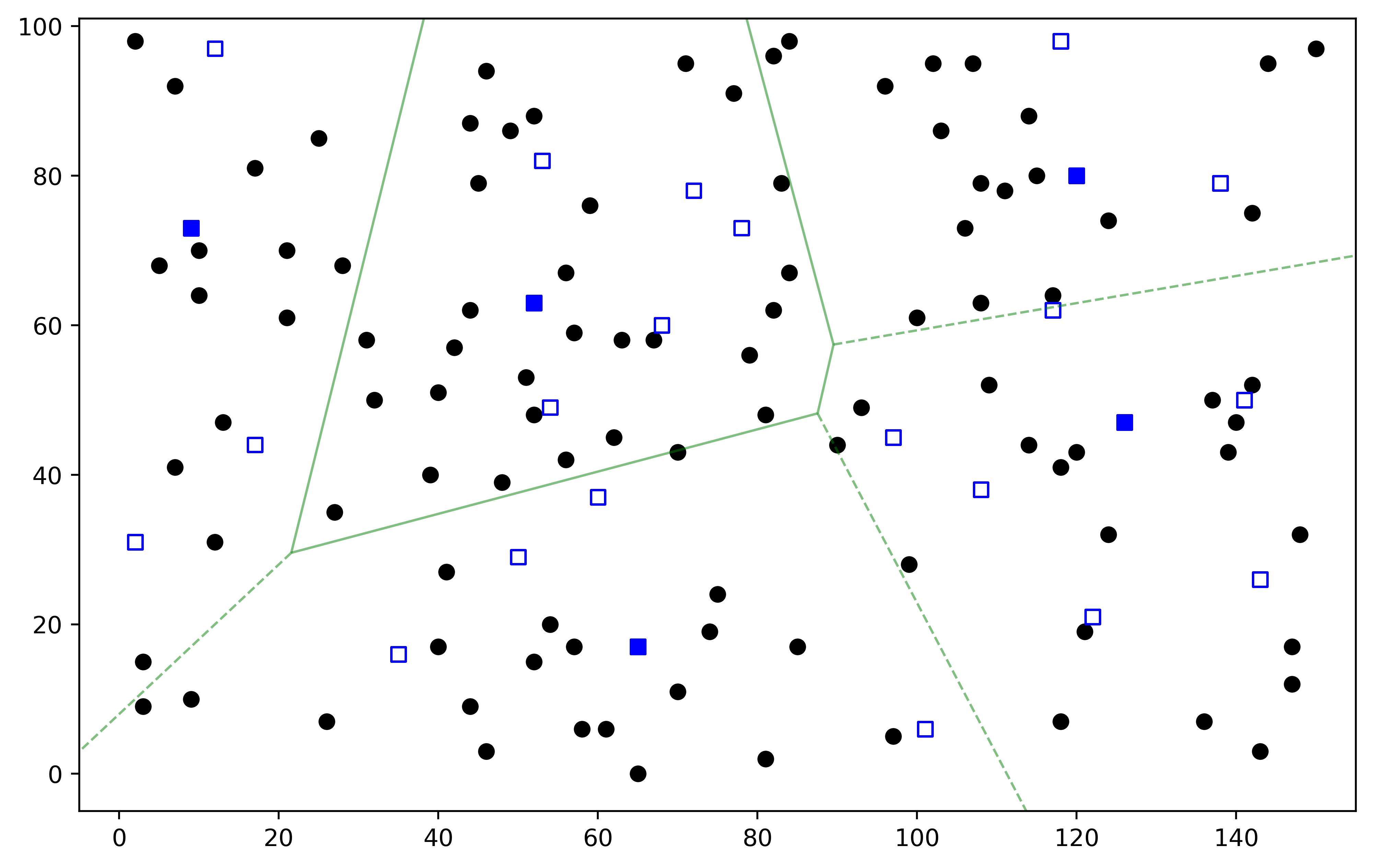}
        \caption{The obtained solution with minimum $F_2$}
    \end{subfigure}

        \vspace{0.2cm} 
    
     \begin{subfigure}{0.49\textwidth}
        \centering
        \includegraphics[width=\linewidth]{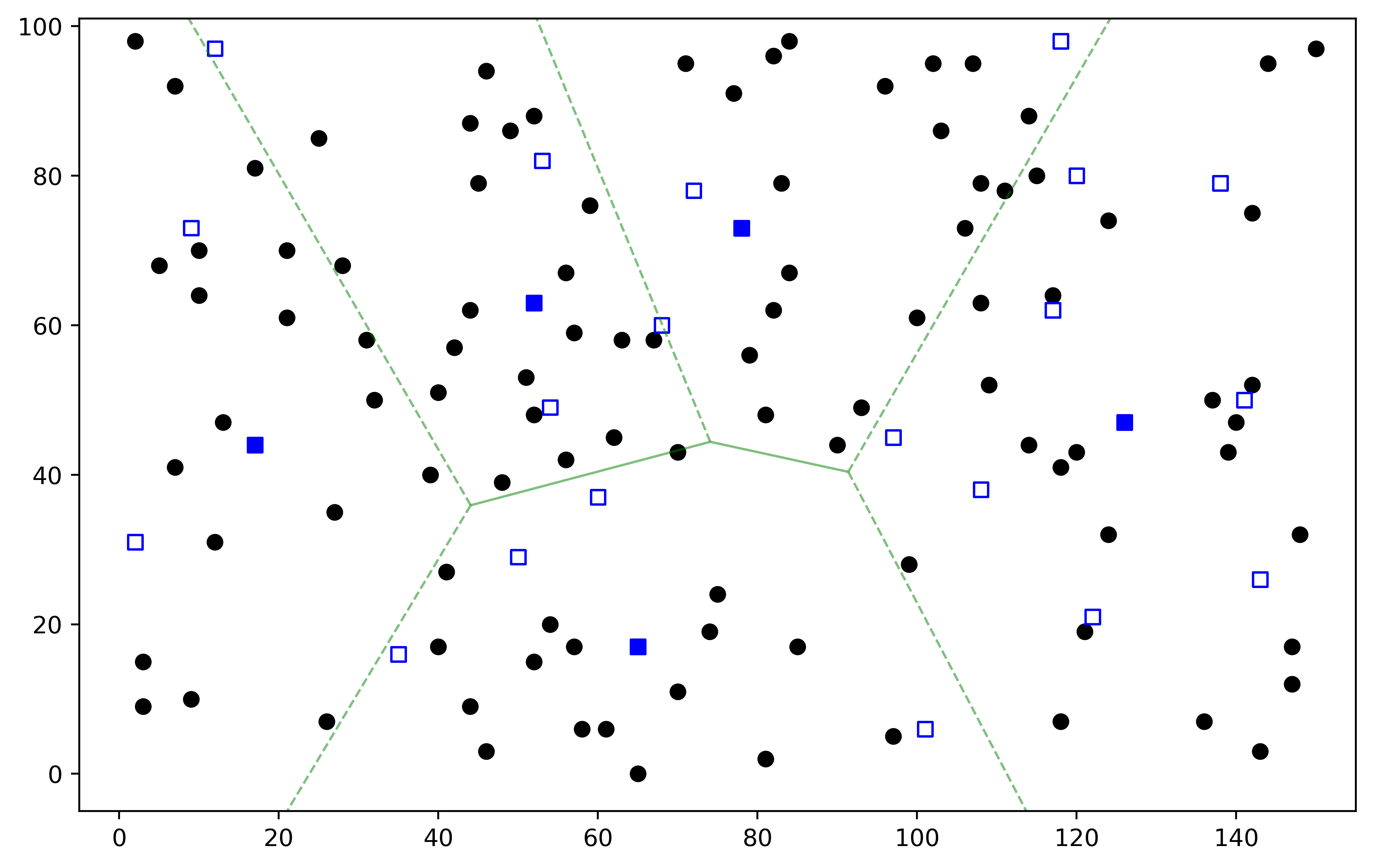}
        \caption{The middle solution among the obtained non-dominated solutions}
    \end{subfigure}
    \hfill
    \begin{subfigure}{0.49\textwidth}
        \centering
        \includegraphics[width=\linewidth]{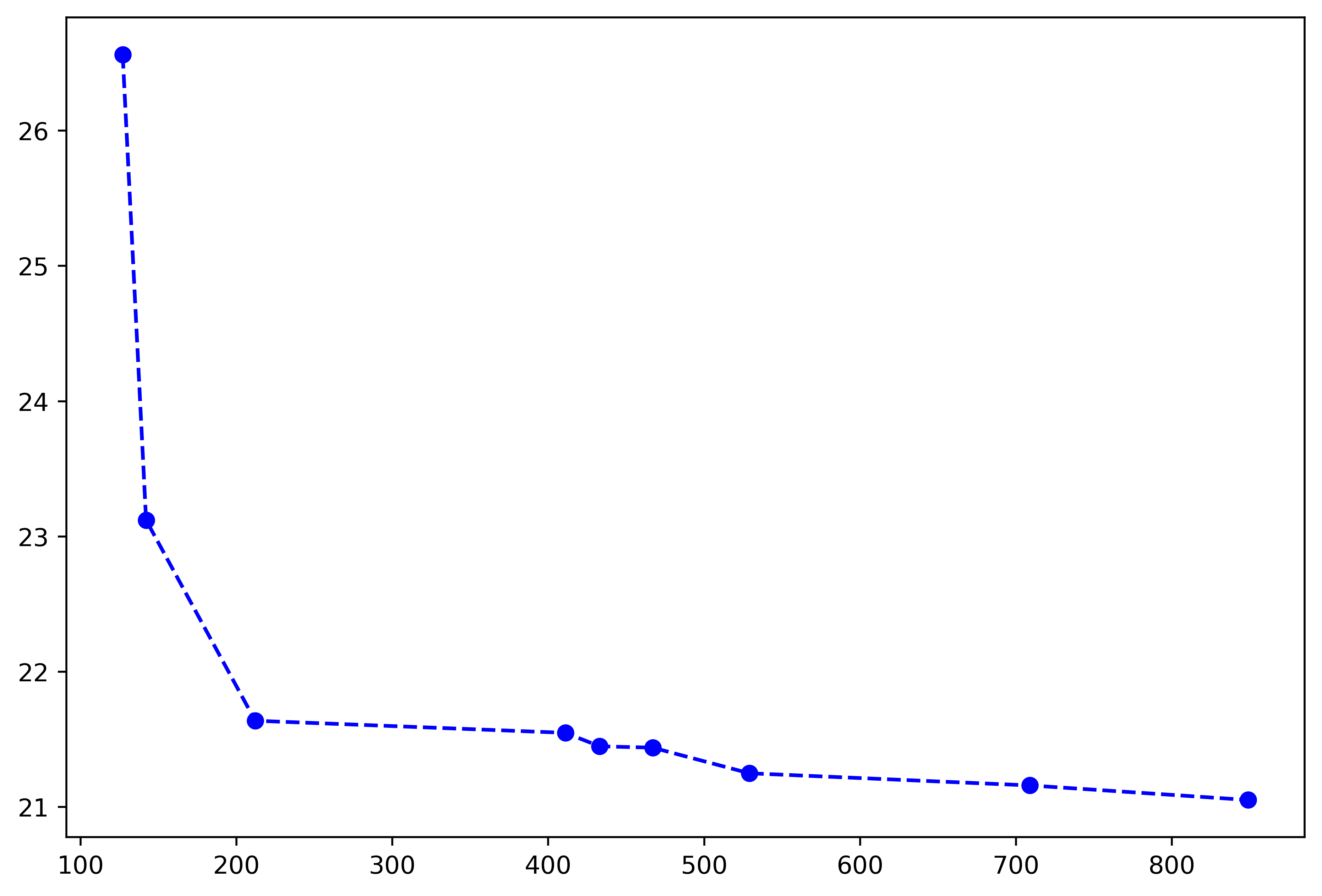}
        \caption{Visualization of all obtained non-dominated solutions in the objective space}
    \end{subfigure}

    \caption{Obtained non-dominated solutions for an instance (100,25,5) by TCLA
    \label{Figk5}}
\end{figure}

\begin{figure}[h]
    \centering
    \begin{subfigure}{0.49\textwidth}
        \centering
        \includegraphics[width=\linewidth]{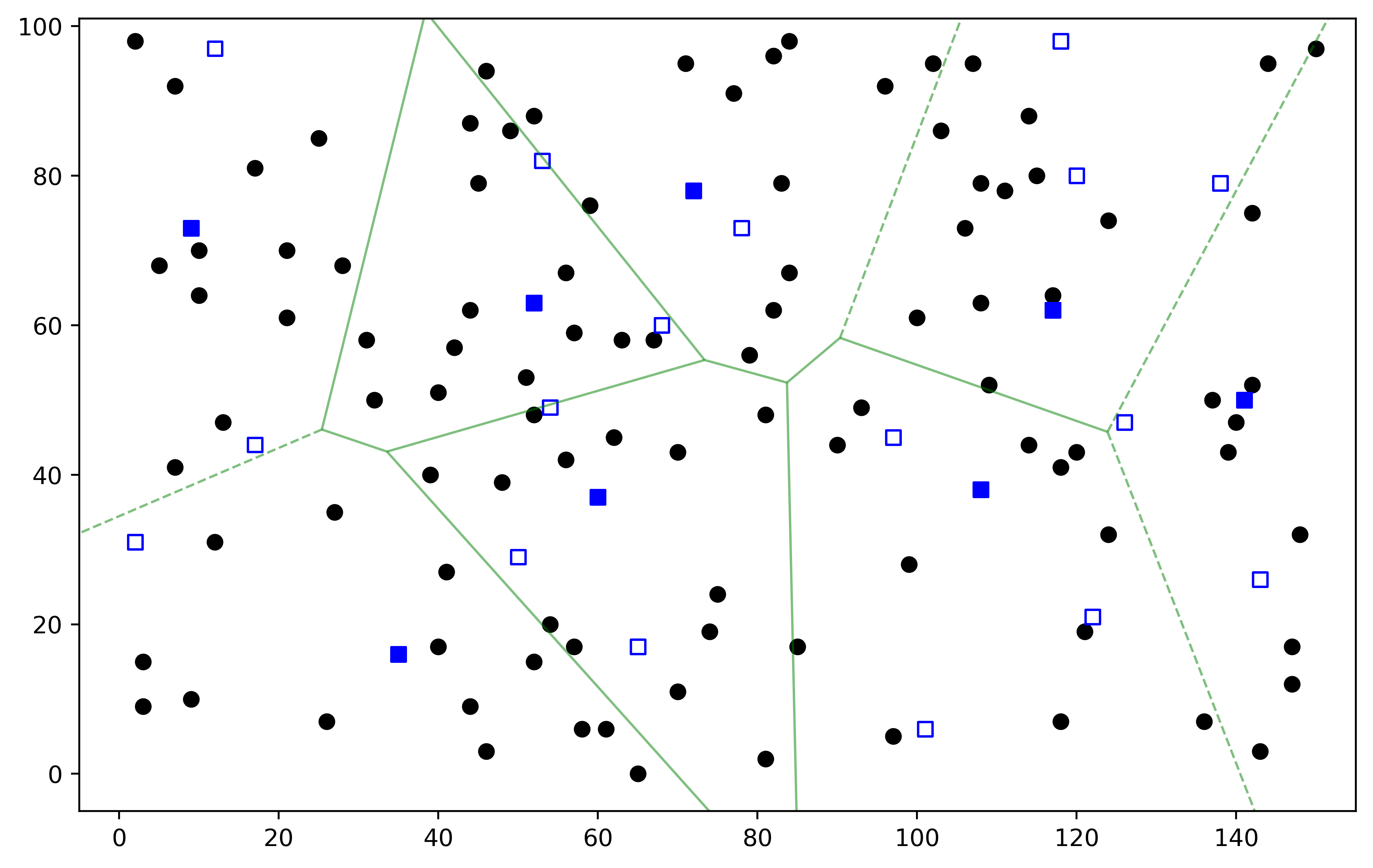}
        \caption{The obtained solution with minimum $F_1$}
    \end{subfigure}
    \hfill
    \begin{subfigure}{0.49\textwidth}
        \centering
        \includegraphics[width=\linewidth]{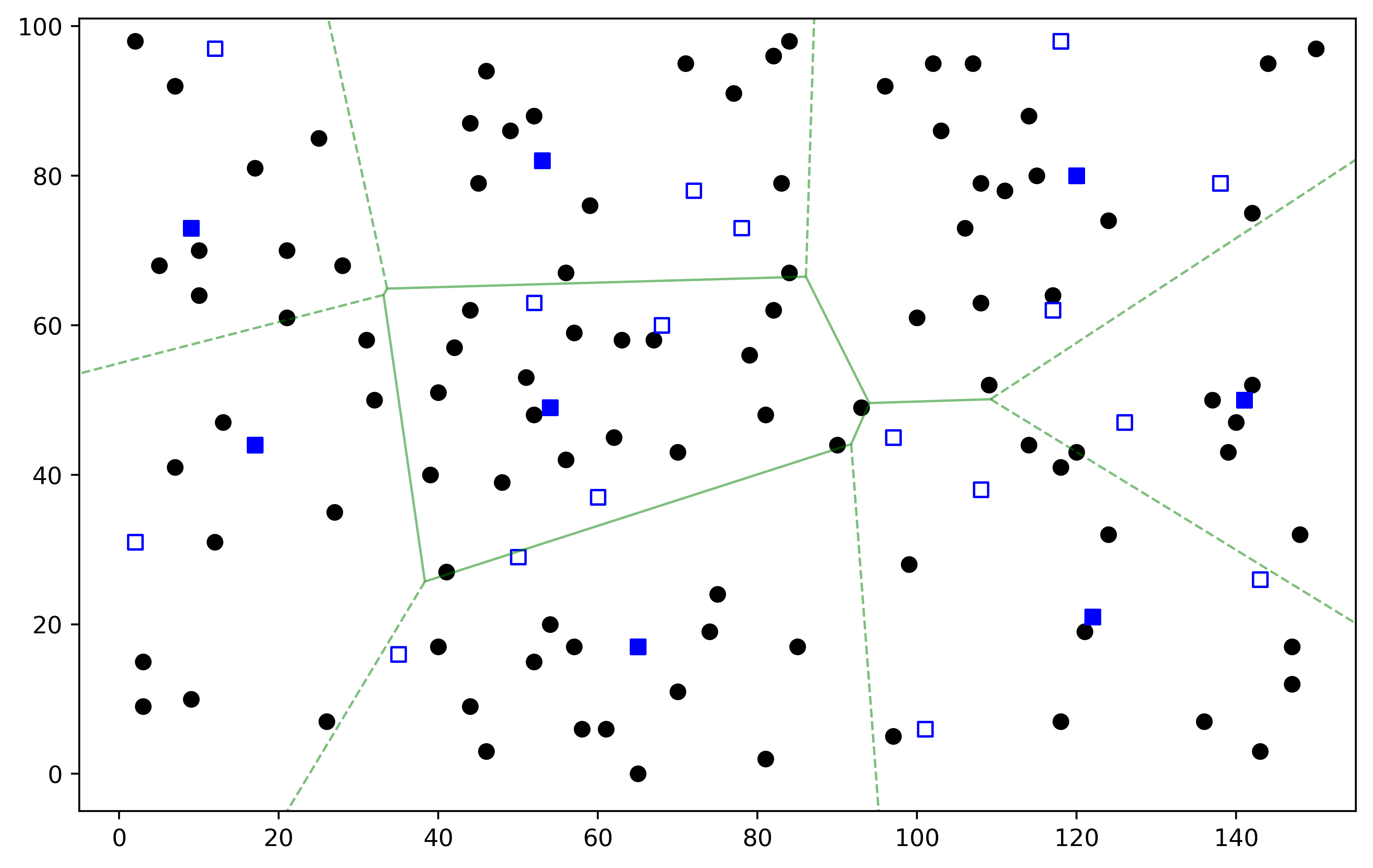}
        \caption{The obtained solution with minimum $F_2$}
    \end{subfigure}

        \vspace{0.2cm} 
    
     \begin{subfigure}{0.49\textwidth}
        \centering
        \includegraphics[width=\linewidth]{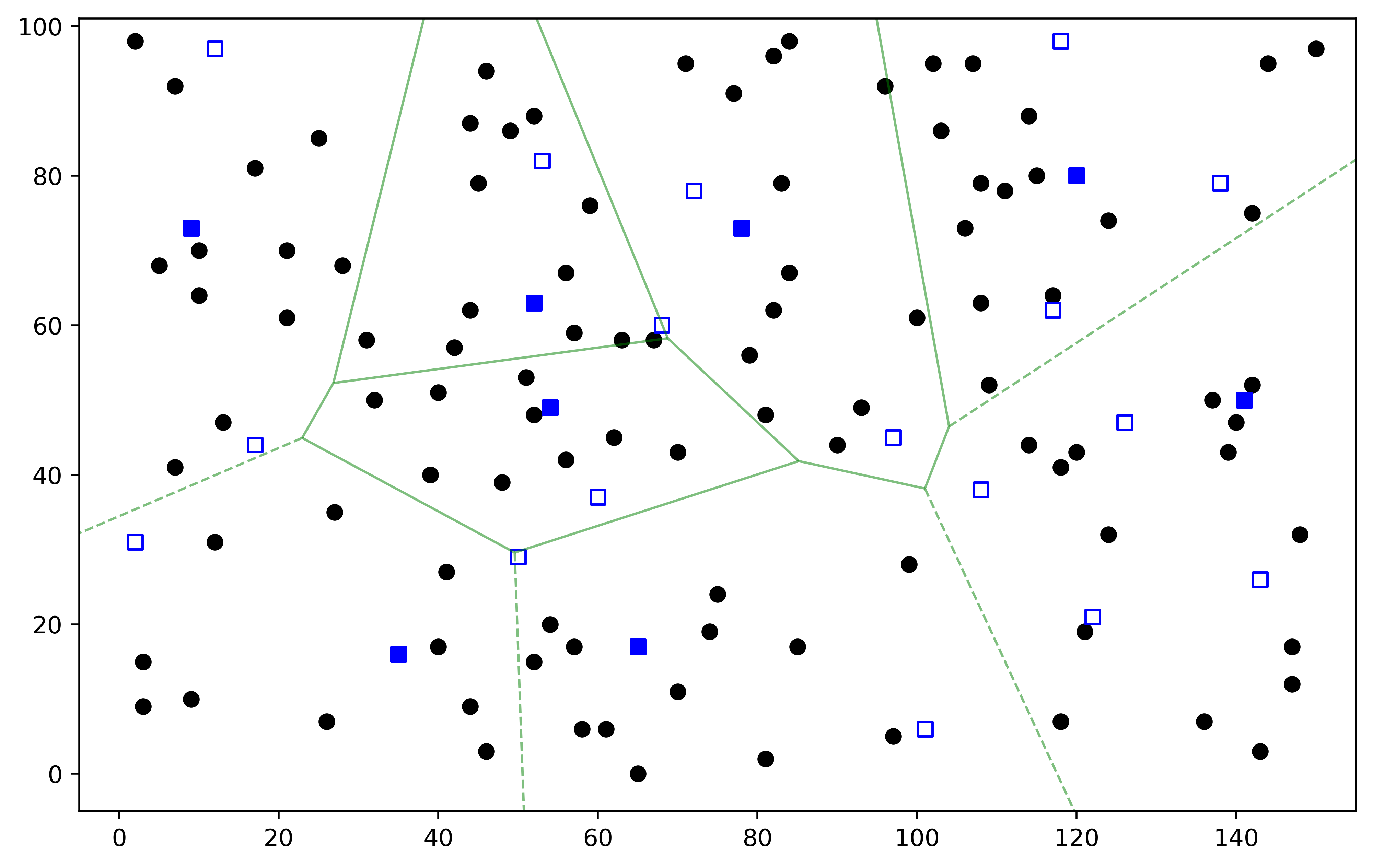}
        \caption{The middle solution among the obtained non-dominated solutions}
    \end{subfigure}
    \hfill
    \begin{subfigure}{0.49\textwidth}
        \centering
        \includegraphics[width=\linewidth]{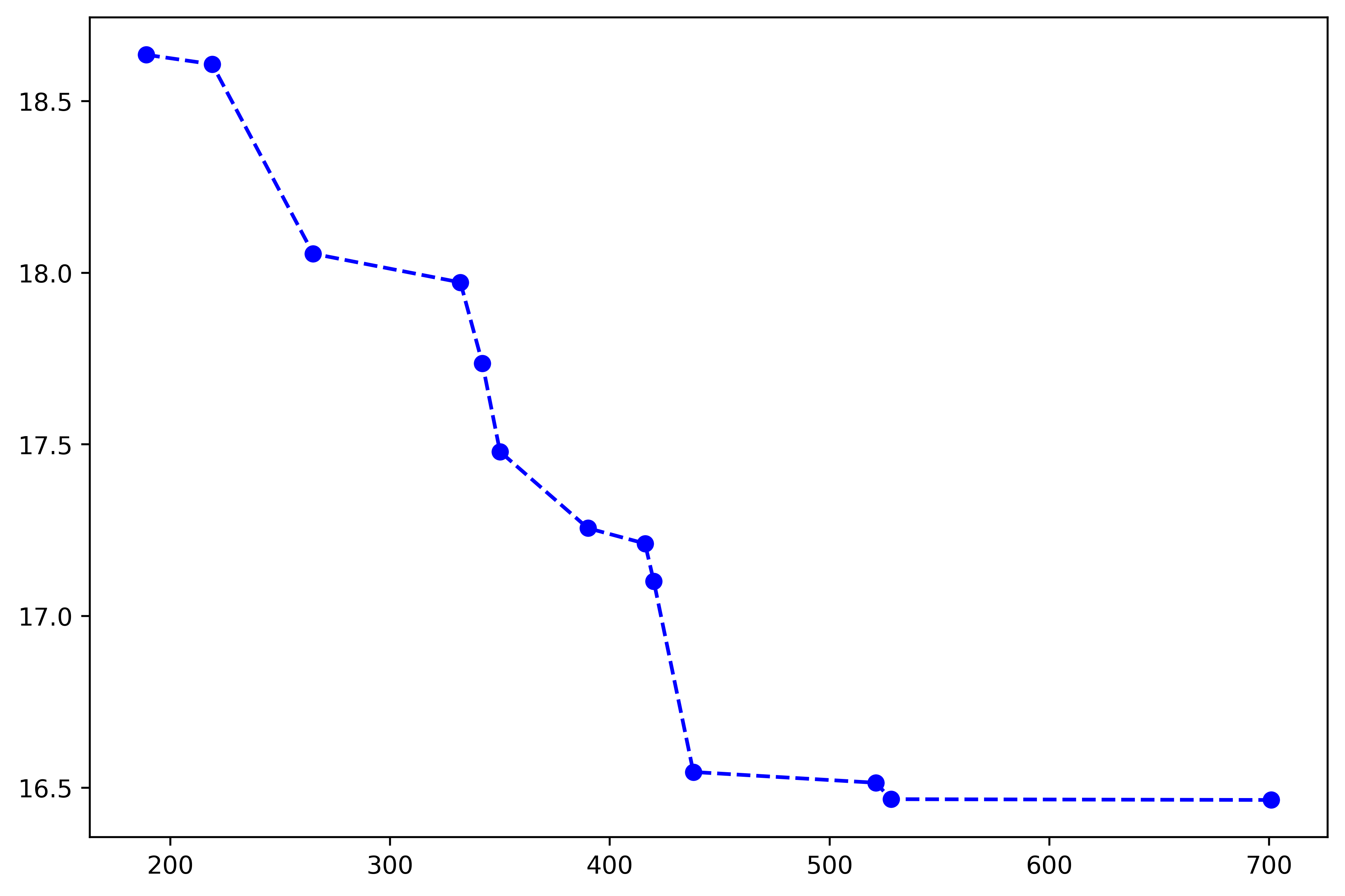}
        \caption{Visualization of all obtained non-dominated solutions in the objective space}
    \end{subfigure}

    \caption{Obtained non-dominated solutions for an instance (100,25,8) by TCLA}
    \label{Figk8}
\end{figure}

\begin{figure}[h]
    \centering
    \begin{subfigure}{0.49\textwidth}
        \centering
        \includegraphics[width=\linewidth]{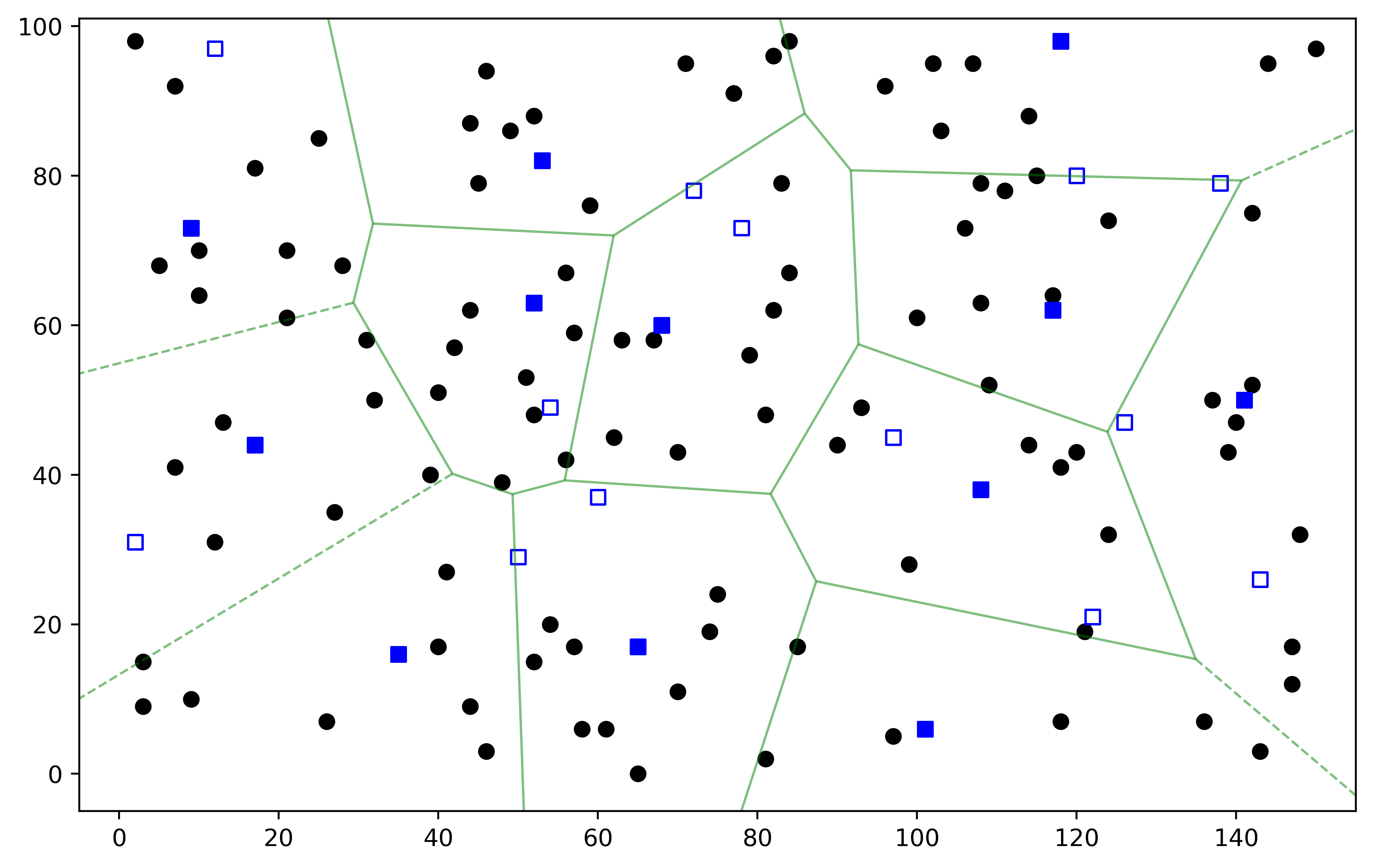}
        \caption{The obtained solution with minimum $F_1$}
    \end{subfigure}
    \hfill
    \begin{subfigure}{0.49\textwidth}
        \centering
        \includegraphics[width=\linewidth]{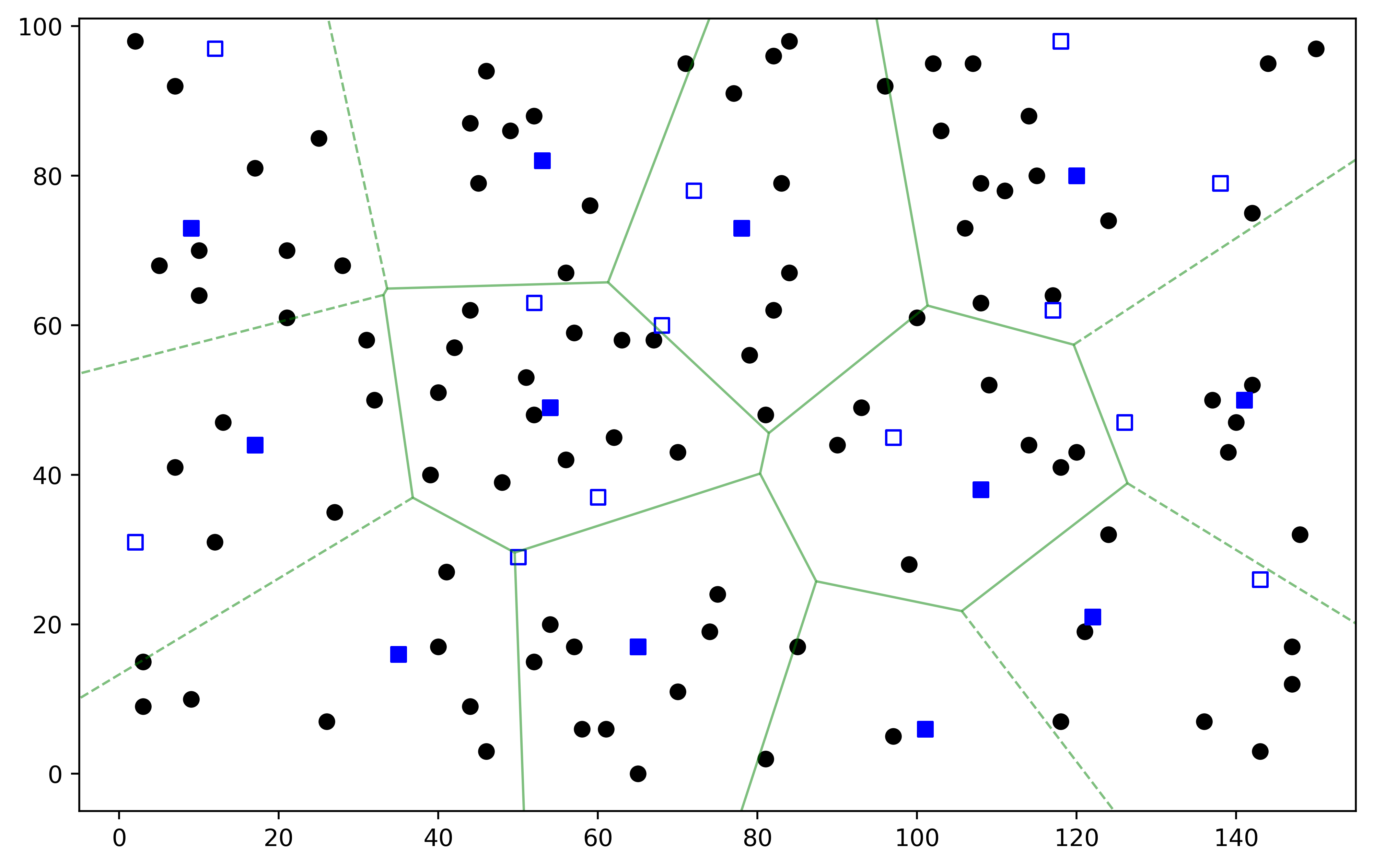}
        \caption{The obtained solution with minimum $F_2$}
    \end{subfigure}

        \vspace{0.2cm} 
    
     \begin{subfigure}{0.49\textwidth}
        \centering
        \includegraphics[width=\linewidth]{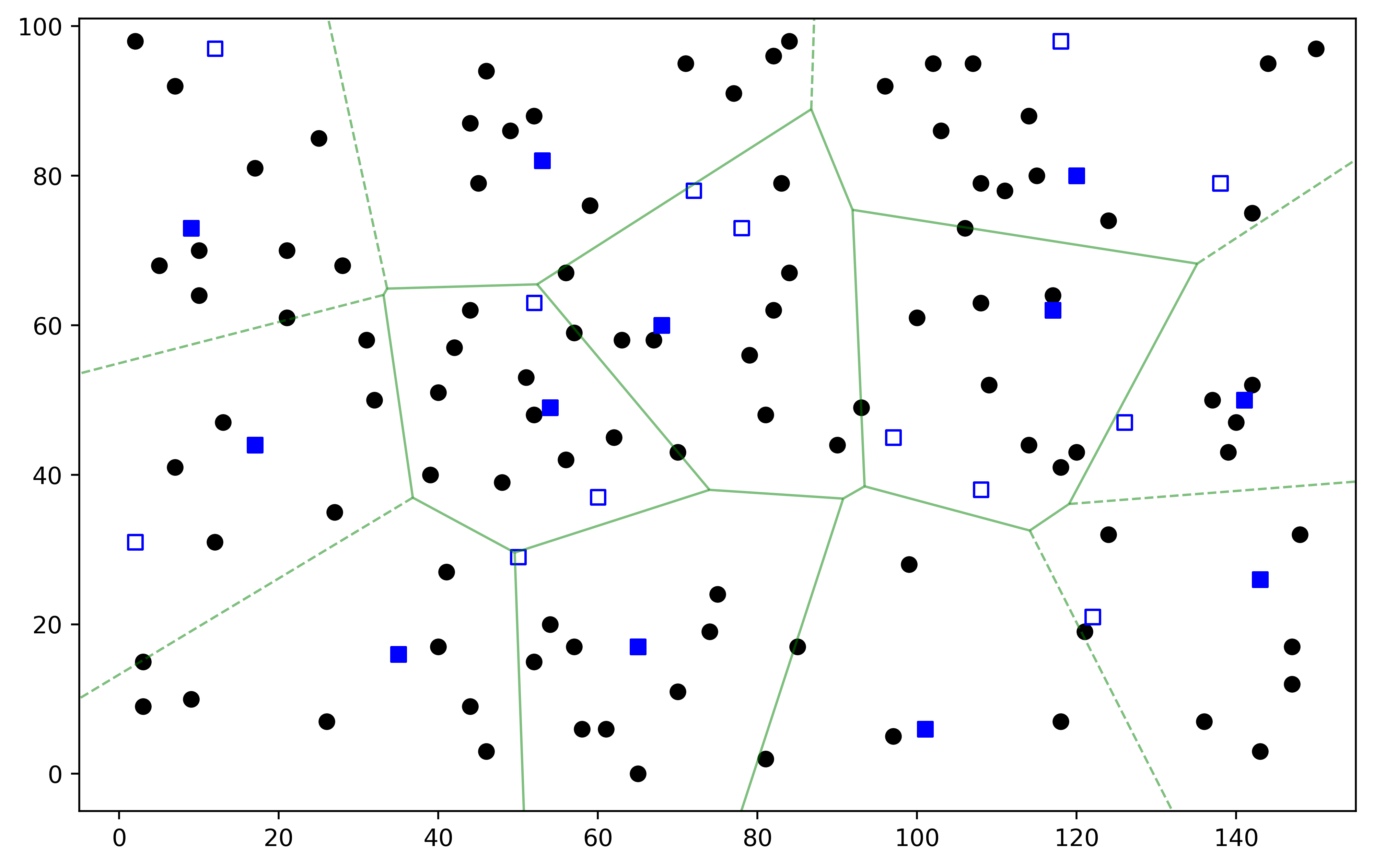}
        \caption{The middle solution among the obtained non-dominated solutions}
    \end{subfigure}
    \hfill
    \begin{subfigure}{0.49\textwidth}
        \centering
        \includegraphics[width=\linewidth]{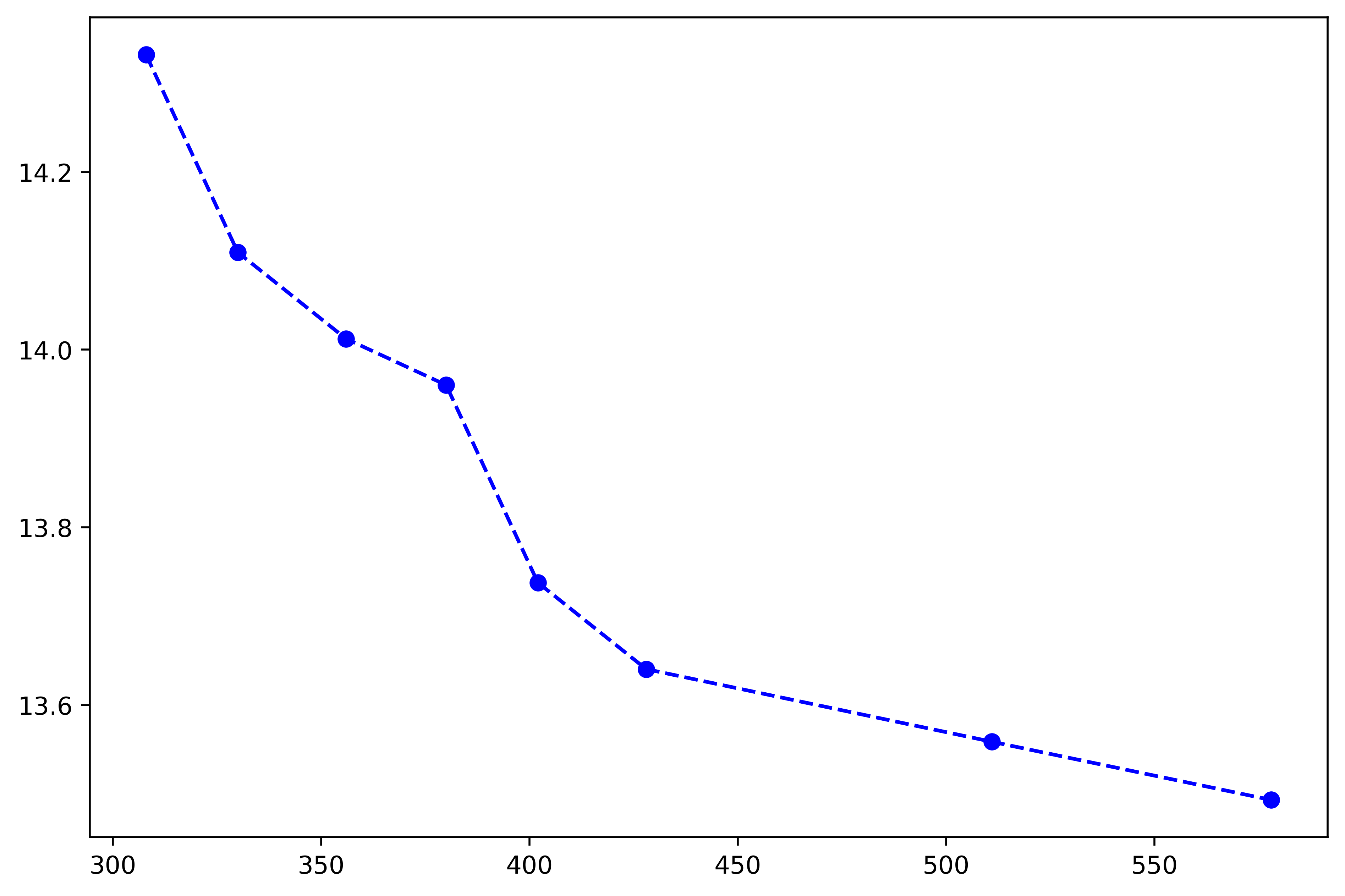}
        \caption{Visualization of all obtained non-dominated solutions in the objective space}
    \end{subfigure}

    \caption{Obtained non-dominated solutions for an instance (100,25,12) by TCLA}
    \label{Figk12}
\end{figure}

\begin{figure}[h]
    \centering
    \begin{subfigure}{0.49\textwidth}
        \centering
        \includegraphics[width=\linewidth]{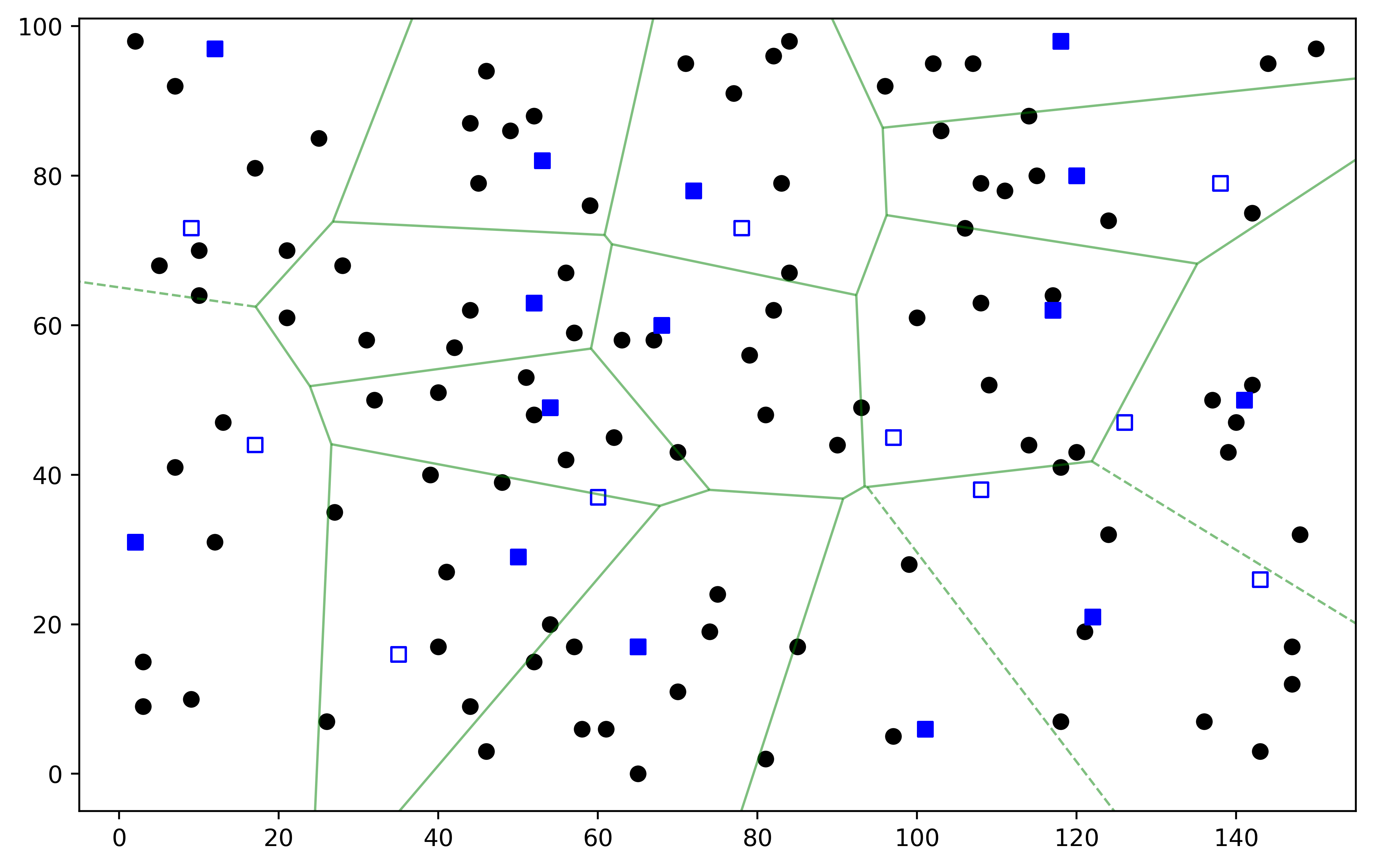}
        \caption{The obtained solution with minimum $F_1$}
    \end{subfigure}
    \hfill
    \begin{subfigure}{0.49\textwidth}
        \centering
        \includegraphics[width=\linewidth]{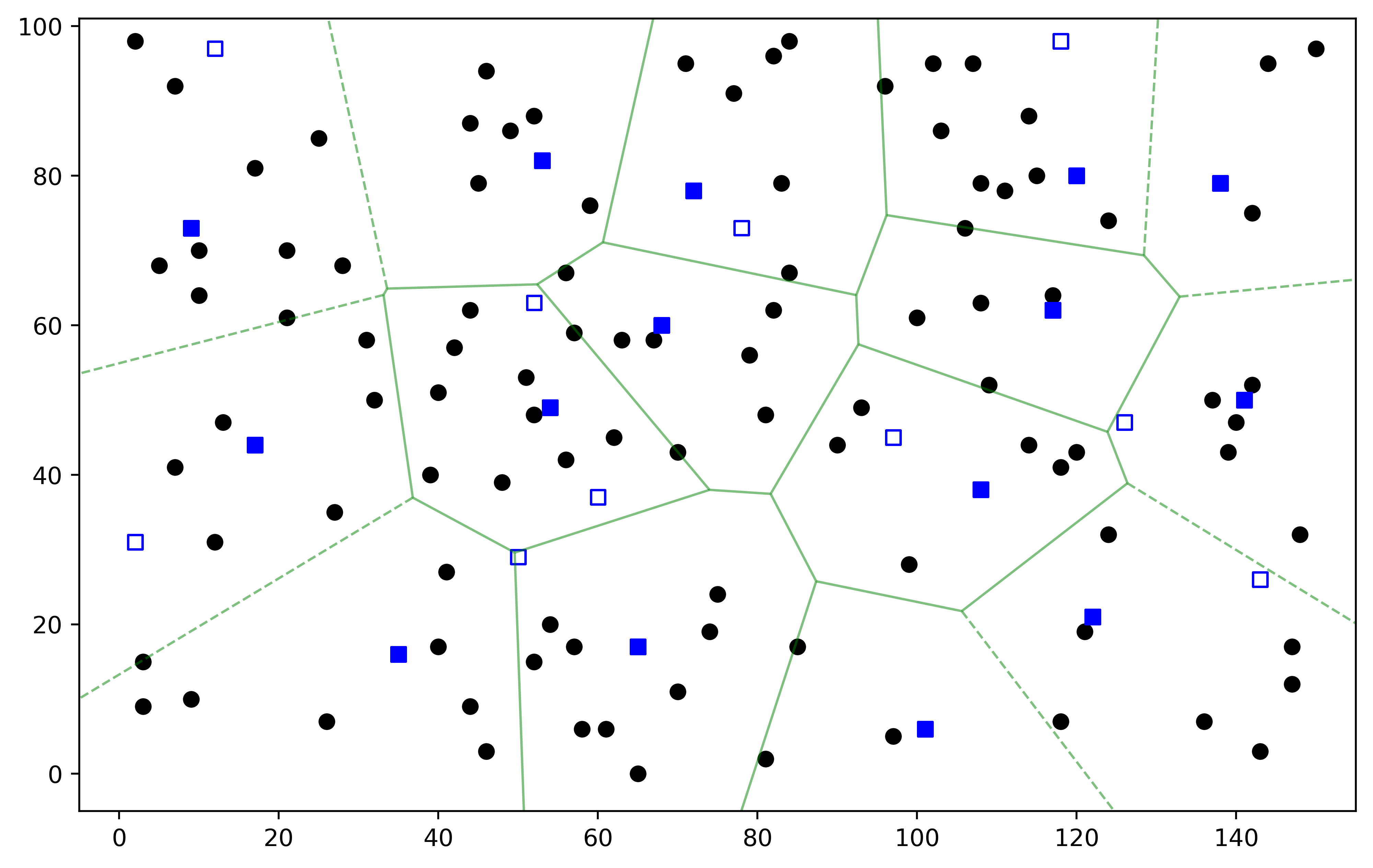}
        \caption{The obtained solution with minimum $F_2$}
    \end{subfigure}

        \vspace{0.2cm} 
    
     \begin{subfigure}{0.49\textwidth}
        \centering
        \includegraphics[width=\linewidth]{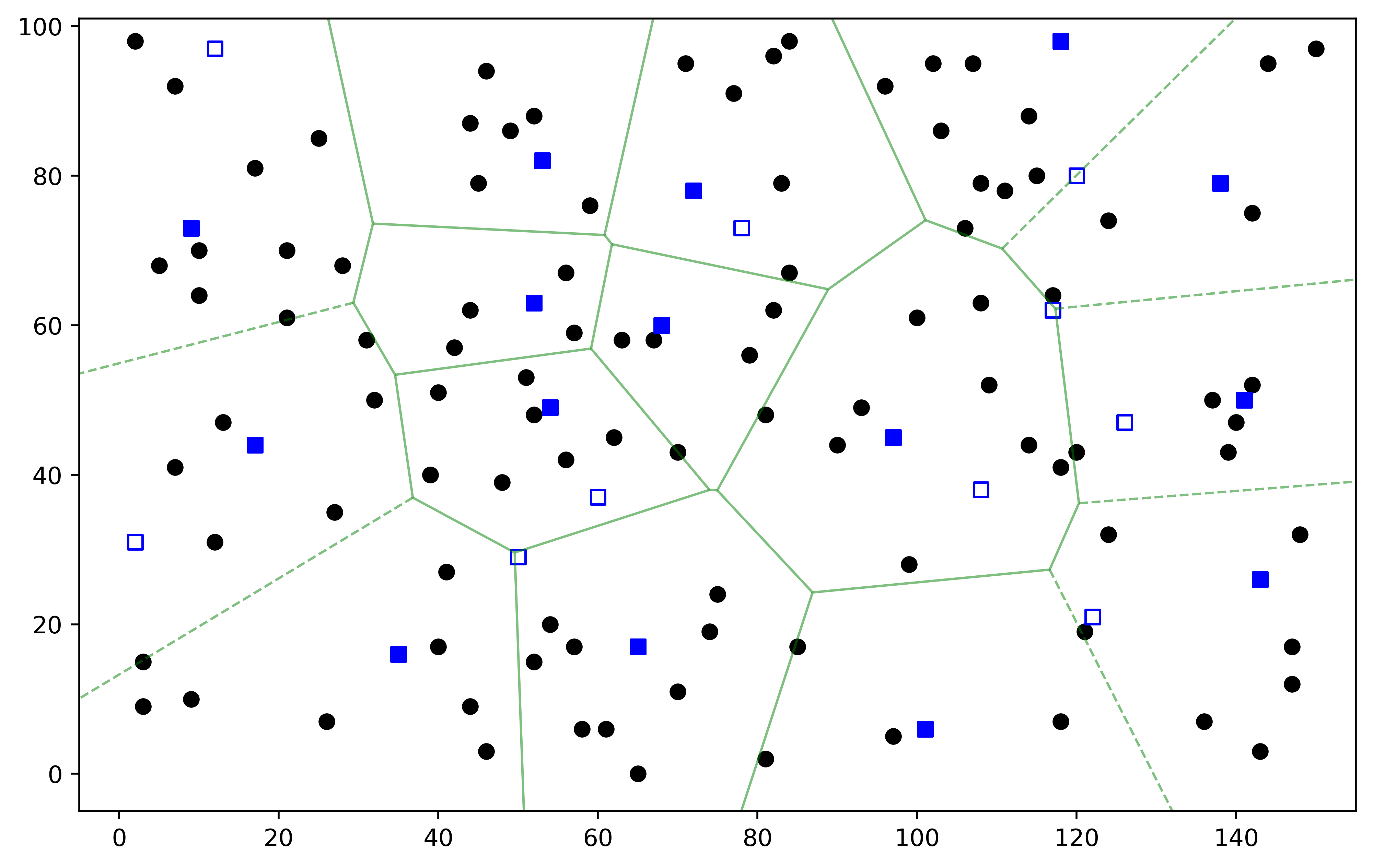}
        \caption{The middle solution among the obtained non-dominated solutions}
    \end{subfigure}
    \hfill
    \begin{subfigure}{0.49\textwidth}
        \centering
        \includegraphics[width=\linewidth]{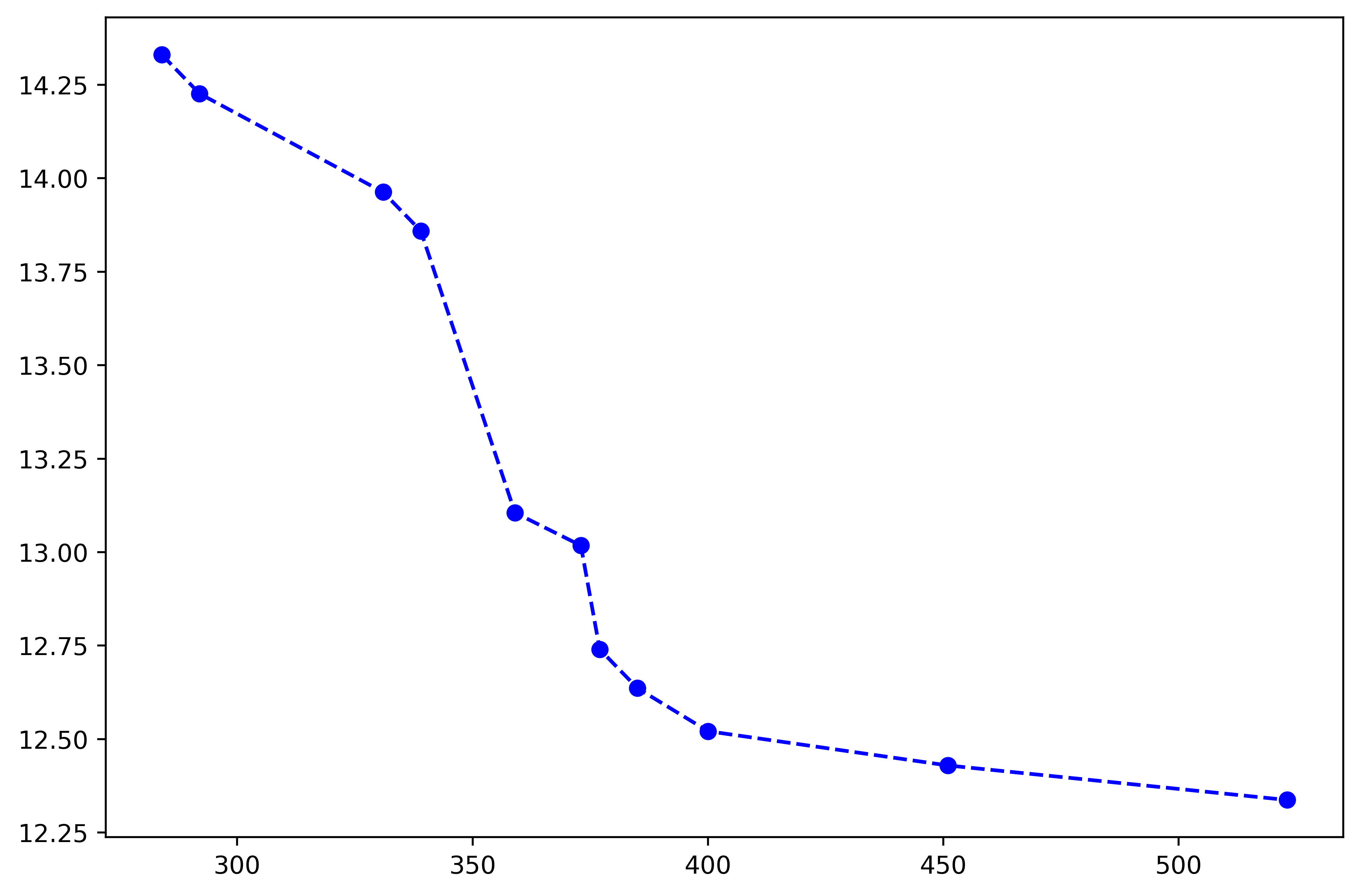}
        \caption{Visualization of all obtained non-dominated solutions in the objective space}
    \end{subfigure}

    \caption{Obtained non-dominated solutions for an instance (100,25,15) by TCLA}
    \label{Figk15}
\end{figure}

Figures \ref{Figk5}--\ref{Figk15} show the results for the instance illustrated in Figure \ref{randenv} for $k=5$, $k=8$, $k=12$ and $k=15$, respectively. In each figure, we select three solutions from the obtained non-dominated set, two extreme solutions and one middle solution. Indeed, we sort the obtained non-dominated solutions according to one of the objectives, i.e., $F_1$, and choose the first (subfigure (\textit{a})), last (subfigure (\textit{b})) and middle (subfigure (\textit{c})) solutions. Also, we depict all the obtained non-dominated solutions in the objective space (subfigure (\textit{a})). The solid blue squares show the selected (opened) test center location in each solution, and for simplicity, we draw the Voronoi edges of them (the green lines). Consequently, the Voronoi region of each selected test center and the demand points that are allocated to each center can be recognized easily. Note that the demand points are weighted (see Figure \ref{randenv}).

The running times of TCLA for $k=5$, $k=8$, $k=12$, and $k=15$ are approximately 7, 24, 41, and 33 seconds, respectively. TCLA finds 9, 13, 8, and 11 non-dominated solutions for $k=5$, $k=8$, $k=12$, and $k=15$, respectively. For $k=5$, the range of $F_1$ values spans from $127$ to $829$, and their corresponding $F_2$ values vary between $23.1$ and $21$. The resulting solution set exhibits a good distribution pattern along the $F_1$ axis. However, there exists a noticeable gap in the $F_2$ values, from $23.1$ to $26.6$, where no solutions are found.
In the cases of $k=8$ and $k=12$, the solution sets exhibit a well-distributed spread in both objective spaces. For $k=8$, the $F_1$ values range from $189$ to $701$, while the $F_2$ values lie in the interval of ($16.5$, $18.6$). Conversely, for $k=12$, the $F_1$ values span from $308$ to $578$, and the $F_2$ values range from $14.3$ to $13.5$.
Lastly, for $k=15$, the obtained solution set covers a range of $F_1$ values from $284$ to $523$ and $F_2$ values between $12.3$ and $14.3$. 

Pareto-optimal solutions play a significant role in aiding decision-makers when selecting an efficient trade-off solution. It is essential to recognize that enhancing one objective often necessitates a trade-off with another objective. The degree of improvement and the associated trade-offs require careful examination. For instance, within the set of non-dominated solutions obtained for the case $k=5$, the third solution with the objective values $F_1 = 212$ and $F_2=21.6$ (refer to Figure \ref{Figk5}-(\textit{d})), stands out as a superior solution, akin to a \textit{knee point}, in comparison to the other solutions within the set.

The parameter $k$ typically stems from budget constraints and the test center's expert limits. Consequently, in addition to comparing sets of Pareto-optimal solutions for a fixed value of $k$, decision-makers can gain insights by observing how the objective values evolve when $k$ is altered. For example, the minimum values of the objective $F_2$, the average traveling distance between the individuals and their closest test center, is improved from 21 to 12.3 when $k$ increases from 5 to 15. 

As we have demonstrated, TCLA successfully identifies a diverse range of non-dominated solutions. However, to comprehensively assess its effectiveness in achieving Pareto optimality, we require comparison results with known Pareto-optimal solutions, which will be discussed in the subsequent section.

\textbf{Comparison results}

We employ the \textit{Set Coverage Metric} (SCM) to assess the Pareto-optimality of the ultimate solutions acquired \cite{zitzler2000comparison}. In the context of two solution sets denoted as $A$ and $B$, the SCM (denoted as $scm(A, B)$) is defined as follows. 

\begin{equation}
scm(A,B) = \frac{|\{\textit{b} \in B | \exists~\textit{a}\in A: a~dominates~\textit{b} \}|}{|B|}.
\label{scm}
\end{equation}
Here, we utilize the notation $|.|$ to represent the cardinality (size) of a set. The metric value $scm(A, B) = 1$ signifies that any solutions within set $B$ are dominated by at least one solution in set $A$. Conversely, when $scm(A, B) = 0$, it implies that no solution in $B$ is dominated by any solutions in $A$. Consequently, when $scm(A, B)$ approaches 1 while $scm(B, A)$ approaches 0, it indicates that solution set $A$ outperforms solution set $B$ in terms of Pareto optimality. In the scenario where the set $A$ comprises Pareto-optimal solutions, it is evident that $scm(B, A) = 0$ holds true for any set $B$, however, $scm(A, B)$ serves as a gauge of set $B$'s effectiveness in achieving Pareto optimality, measuring its efficiency in this regard.
In this part, we employ the Gurobi 5.6.3 optimization solver \cite{gurobi} with a specified parameter of $MIPGap = 1e-3$. Our goal is to identify a collection of solutions that are either Pareto-optimal or very close to being Pareto-optimal. These solutions will subsequently be used as one of the sets in the calculation of $scm(A, B)$.
Similarly, we employ TCLA and identify the resulting non-dominated solutions, which will serve as the other set in the calculation of $scm(A, B)$. It is important to note that the solutions obtained through Gurobi may not necessarily be Pareto-optimal due to the presence of an optimality gap. This gap signifies the difference between the best-known integer solution and the current best solution discovered during Gurobi's branch-and-bound process. The optimality gap is expressed as a percentage of the objective function value and plays a crucial role in balancing the trade-off between solution quality and solution computation time.

The metric $scm(A, B)$ calculates the number of solutions within a set $B$ that are dominated by at least one solution in set $A$. However, it lacks the ability to quantify "to what extent" a solution $b \in B$ may be dominated by a solution $a \in A$. To address this concern, we extend the concept of the \textit{approximation factor} from single-objective optimization to bi-objective optimization as follows: If we represent the objective values of a solution $b$ from the set $B$ as $F_1(b)$ and $F_2(b)$, when solution $a$ dominates solution $b$, we have

\begin{equation}
[ F_1(a) < F_1(b)~~and~~ F_2(a) \leq F_2(b) ]
~~~or~~~
[ F_1(a) \leq F_1(b)~~and~~ F_2(a) < F_2(b) ]
\end{equation}

Now, let $\alpha$ and $\beta$ be the smallest values that satisfy the following equations.

\begin{equation}
F_1(a) \geq (1-\alpha) F_1(b)~~and~~ F_2(a) \geq (1-\beta) F_2(b)
\label{augment}
\end{equation}

Indeed, parameter $\alpha$ ($\beta$) signifies the percentage by which solution $b$ must enhance its performance to surpass solution $a$ in objective $F_1$ ($F_2$). Consequently, decreasing $\alpha$ percent in the direction of $F_1$, or $\beta$ percent in the direction of $F_2$, will render solution $b$ no longer dominated by solution $a$. Furthermore, augmenting in both directions simultaneously will yield a stronger solution. Now, identifying the minimum pair of $(\alpha, \beta)$ that satisfies equation Eq. (\ref{augment}) for all solutions $a \in A$ reflects the quality of the solution $b$. Furthermore, extending this fact to all solutions $b \in B$ and choosing the maximum $\alpha$ and $\beta$ values among them will represent the quality of solution set $B$ in comparison to solution set $A$. Let us denote this metric by $\alpha\beta(A,B)$.


To utilize the Gurobi solver, we implement the $\epsilon$-constraint methodology, as elucidated in Section \ref{epsiloncons}.
We establish a range encompassing the lower and upper bounds for the feasible values of $F_2$, denoting the average travel distance between individuals and their nearest open center. Subsequently, based on the desired quantity of Pareto-optimal solutions, we evenly select various $\epsilon$ values from this specified interval. 
For example, for an interval $[left,right]$, if we are interested in finding at most $h$ Pareto-optimal solutions, we run the Gurobi solver for $\epsilon= left + i(\frac{right-left}{h-1})$, for $i=0,1,2,...,h-1$.
It's worth noting that for certain $\epsilon$ values, particularly those close to the lower boundary of the interval, no feasible solutions may be attainable. Additionally, some of these $\epsilon$ values may yield identical optimal solutions, indicating that in the bi-objective space, the solutions obtained from larger $\epsilon$ values are dominated by those obtained from smaller ones. Consequently, we run the Gurobi solver for all $h$ sampled $\epsilon$ values independently and subsequently report solely the Pareto-optimal solutions that have been ultimately obtained.

We generate random instances of varying sizes for the test center location problem and, for each instance, employ both TCLA and Gurobi solver using the described $\epsilon$-constraint approach. For each run, we provide the number of non-dominated solutions obtained by TCLA, the number of Pareto optimal solutions obtained by Gurobi, and the set coverage metric $scm(A,B)$, where $A$ and $B$ are the obtained solutions by Gurobi and TCLA, respectively. In addition, we report the running time (in seconds) for both approaches. It's worth noting that TCLA can find a non-dominated set in a single run, whereas Gurobi executes separately for each $\epsilon$ value, resulting in one Pareto-optimal solution per run. Therefore, for Gurobi runs, we report two types of running times: the total running time for finding all Pareto-optimal solutions, and the average running time to discover a single Pareto-optimal solution, excluding runs that yield no feasible solutions.

\begin{figure}
    \centering
    \begin{subfigure}{0.24\textwidth}
        \centering
        \includegraphics[width=\linewidth]{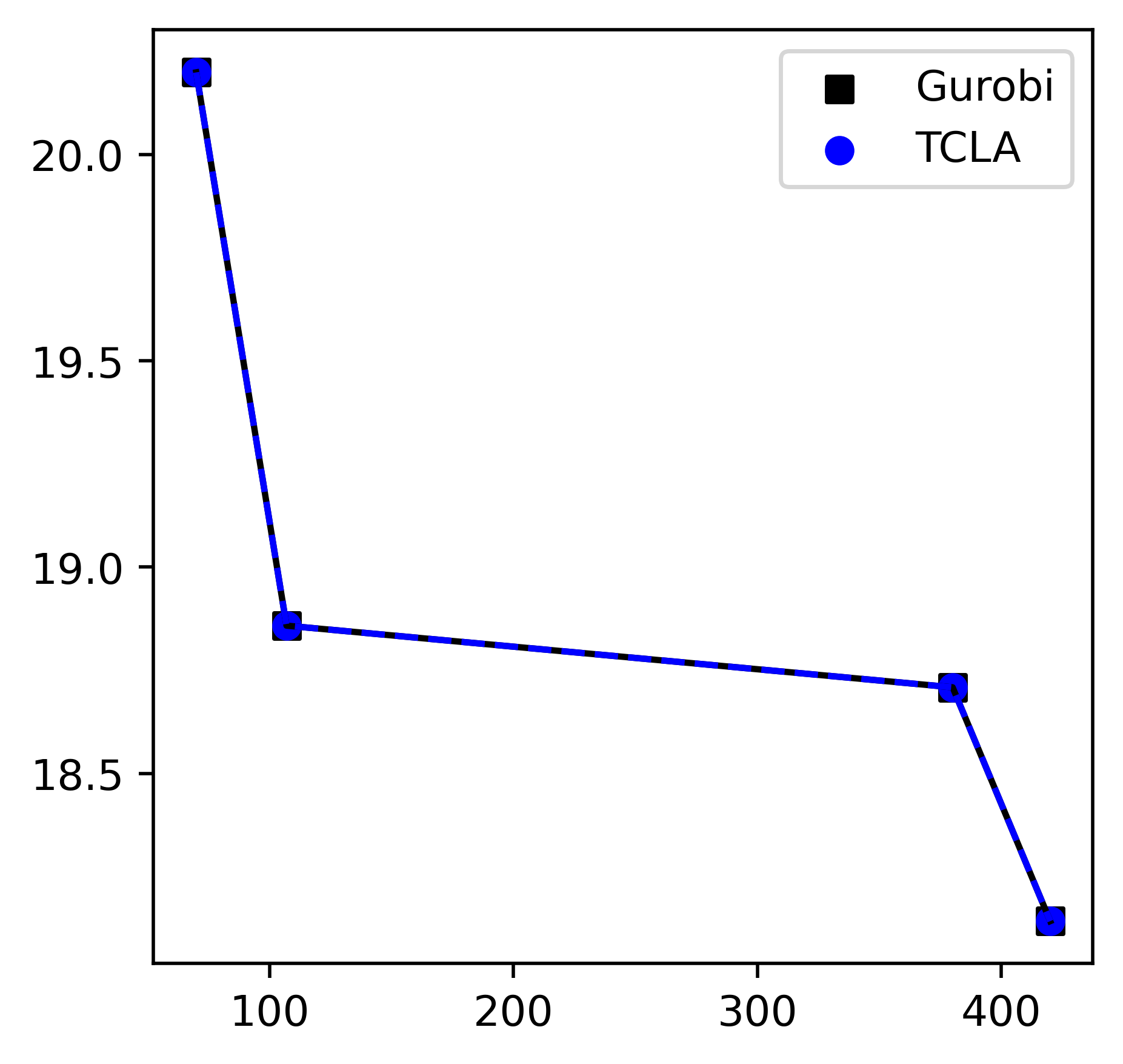}
        \caption*{instance (40,20,5)}
    \end{subfigure}
    \begin{subfigure}{0.24\textwidth}
        \centering
        \includegraphics[width=\linewidth]{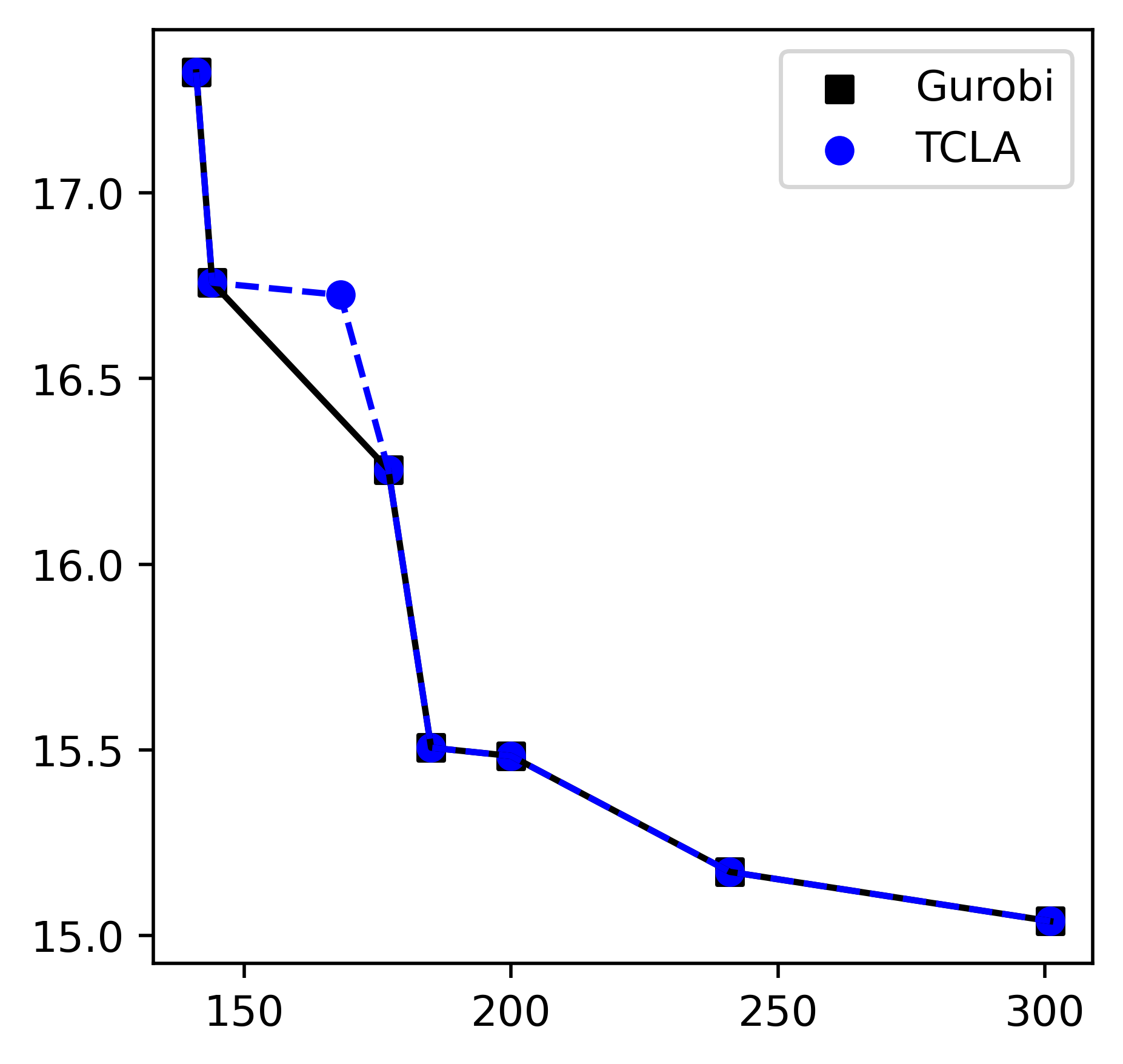}
        \caption*{instance (40,20,8)}
    \end{subfigure}
     \begin{subfigure}{0.24\textwidth}
        \centering
        \includegraphics[width=\linewidth]{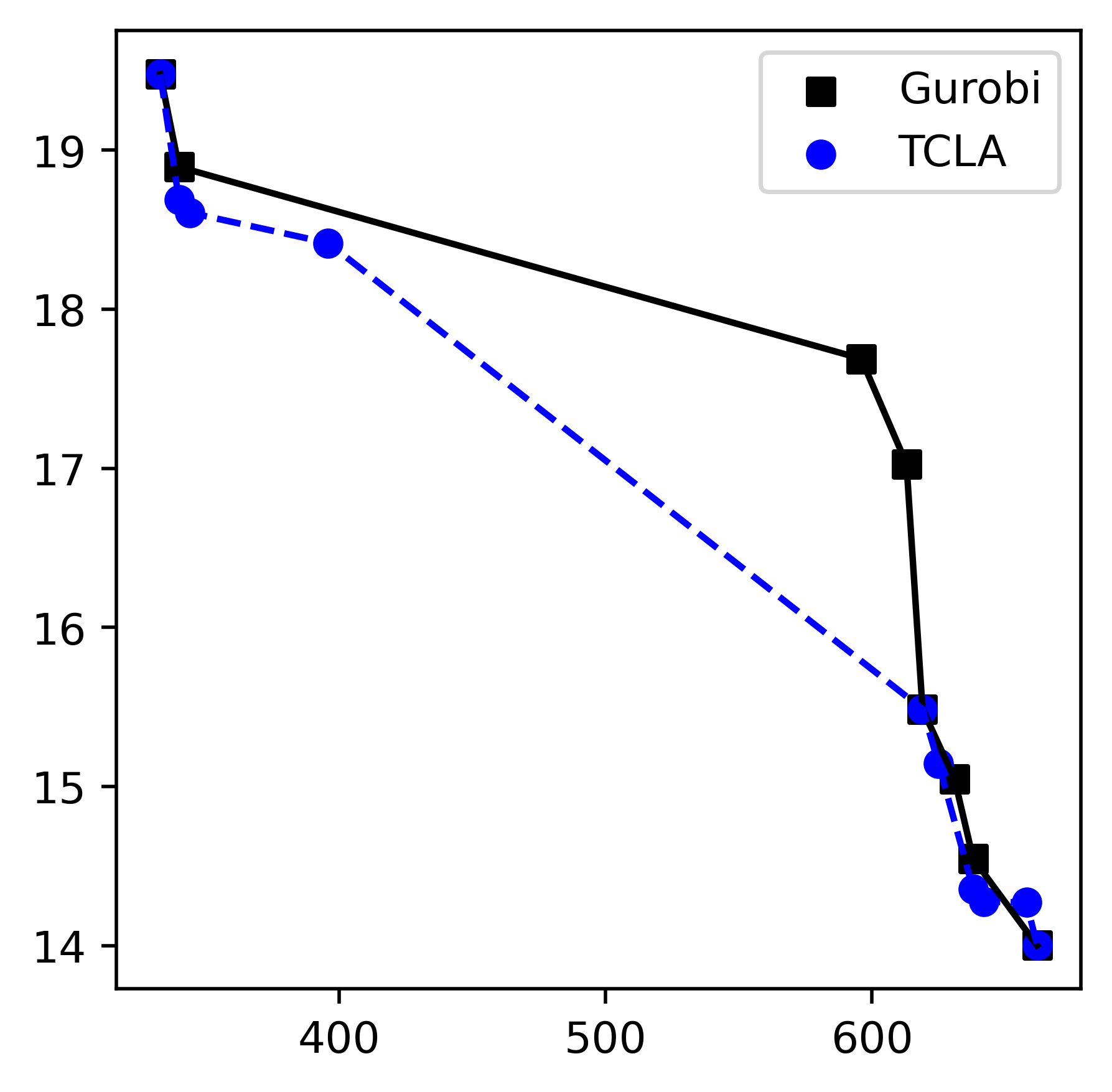}
        \caption*{instance (40,20,10)}
    \end{subfigure}
    \begin{subfigure}{0.24\textwidth}
        \centering
        \includegraphics[width=\linewidth]{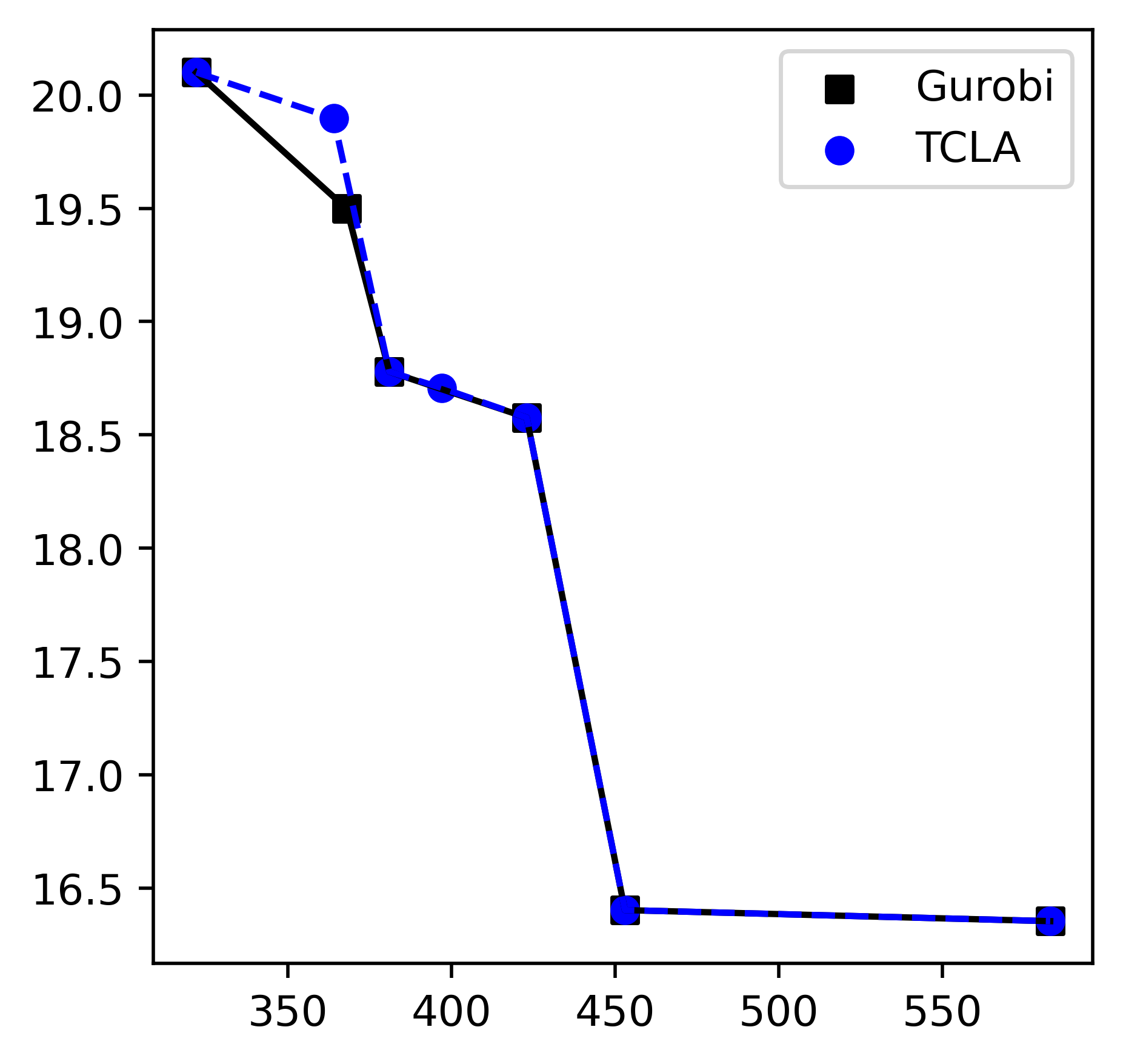}
        \caption*{instance (40,20,12)}
    \end{subfigure}

    \vspace{0.2cm} 

   \begin{subfigure}{0.24\textwidth}
        \centering
        \includegraphics[width=\linewidth]{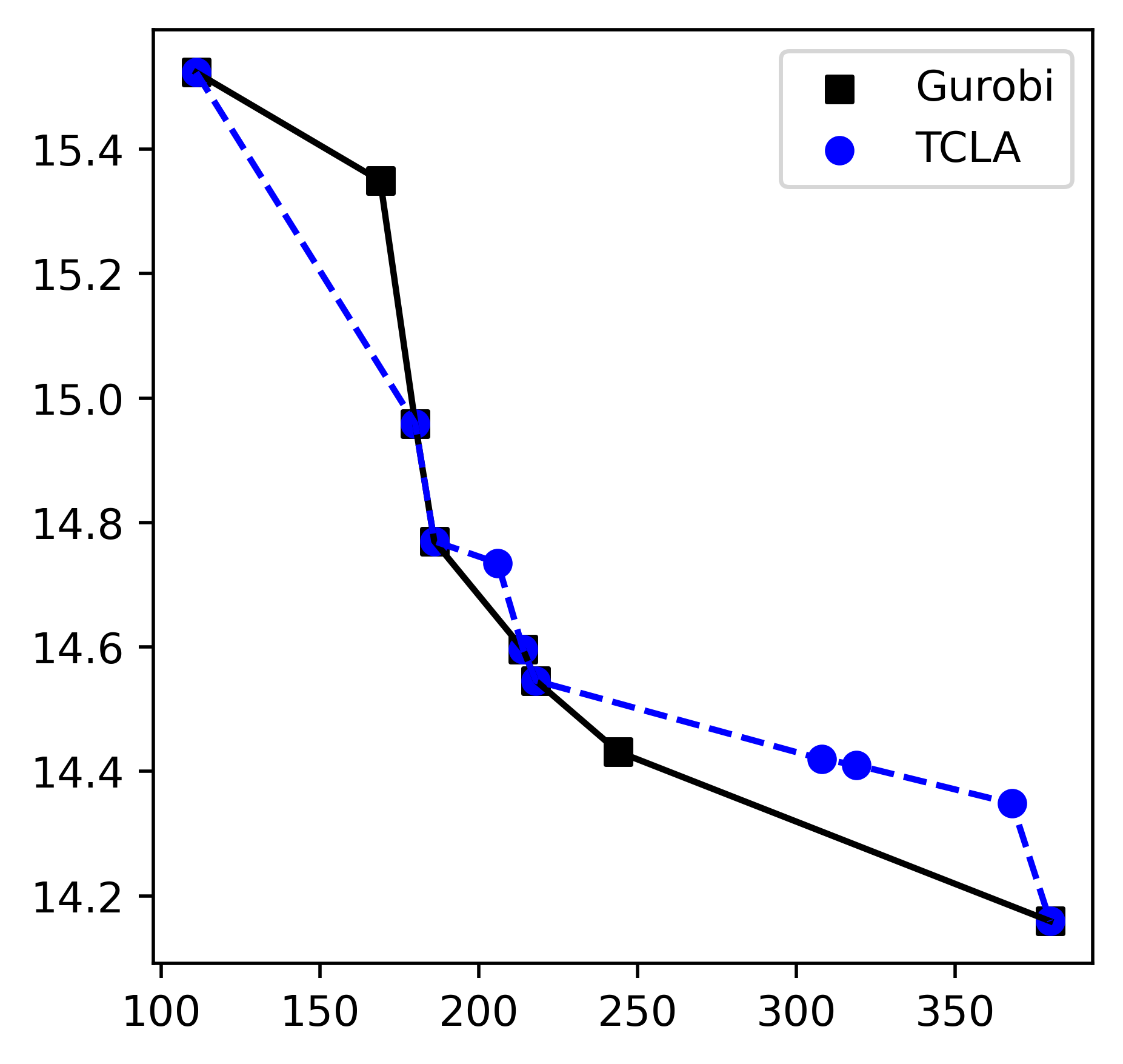}
        \caption*{instance (100,40,10)}
    \end{subfigure}
    \begin{subfigure}{0.24\textwidth}
        \centering
        \includegraphics[width=\linewidth]{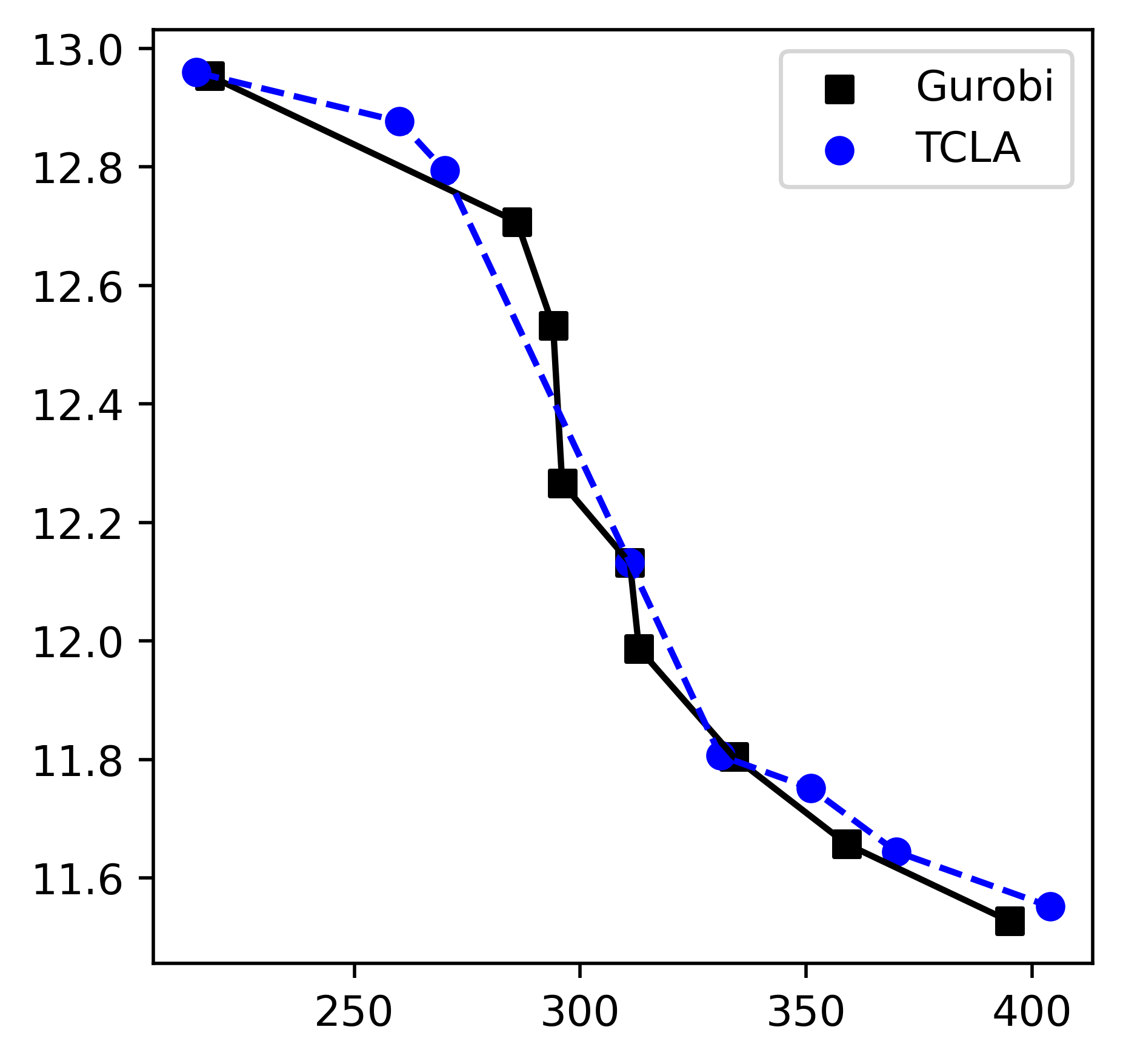}
        \caption*{instance (100,40,15)}
    \end{subfigure}
     \begin{subfigure}{0.24\textwidth}
        \centering
        \includegraphics[width=\linewidth]{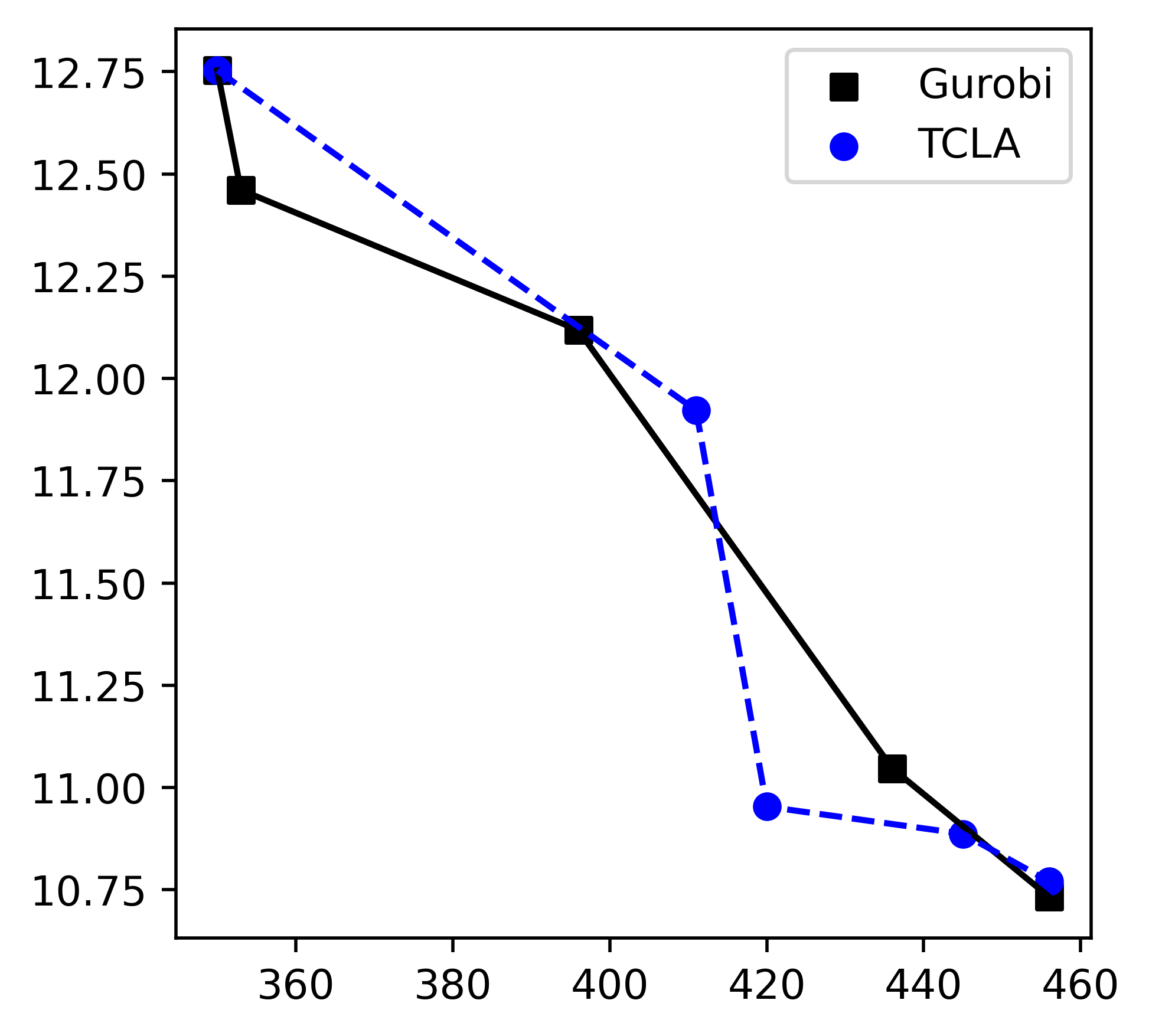}
        \caption*{instance (100,40,20)}
    \end{subfigure}
    \begin{subfigure}{0.24\textwidth}
        \centering
        \includegraphics[width=\linewidth]{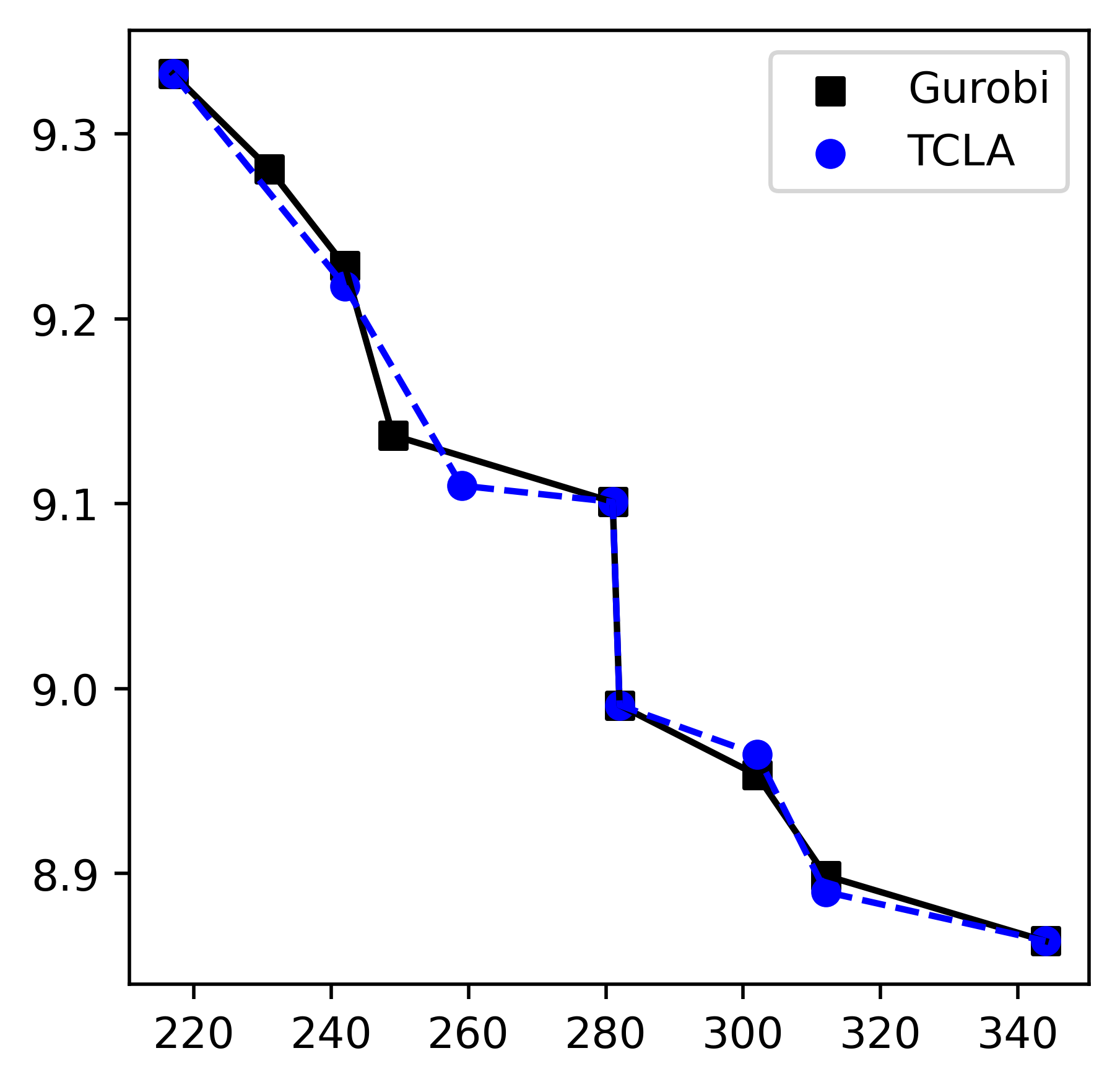}
        \caption*{instance (100,40,25)}
    \end{subfigure}

    \vspace{0.2cm} 

    \begin{subfigure}{0.24\textwidth}
        \centering
        \includegraphics[width=\linewidth]{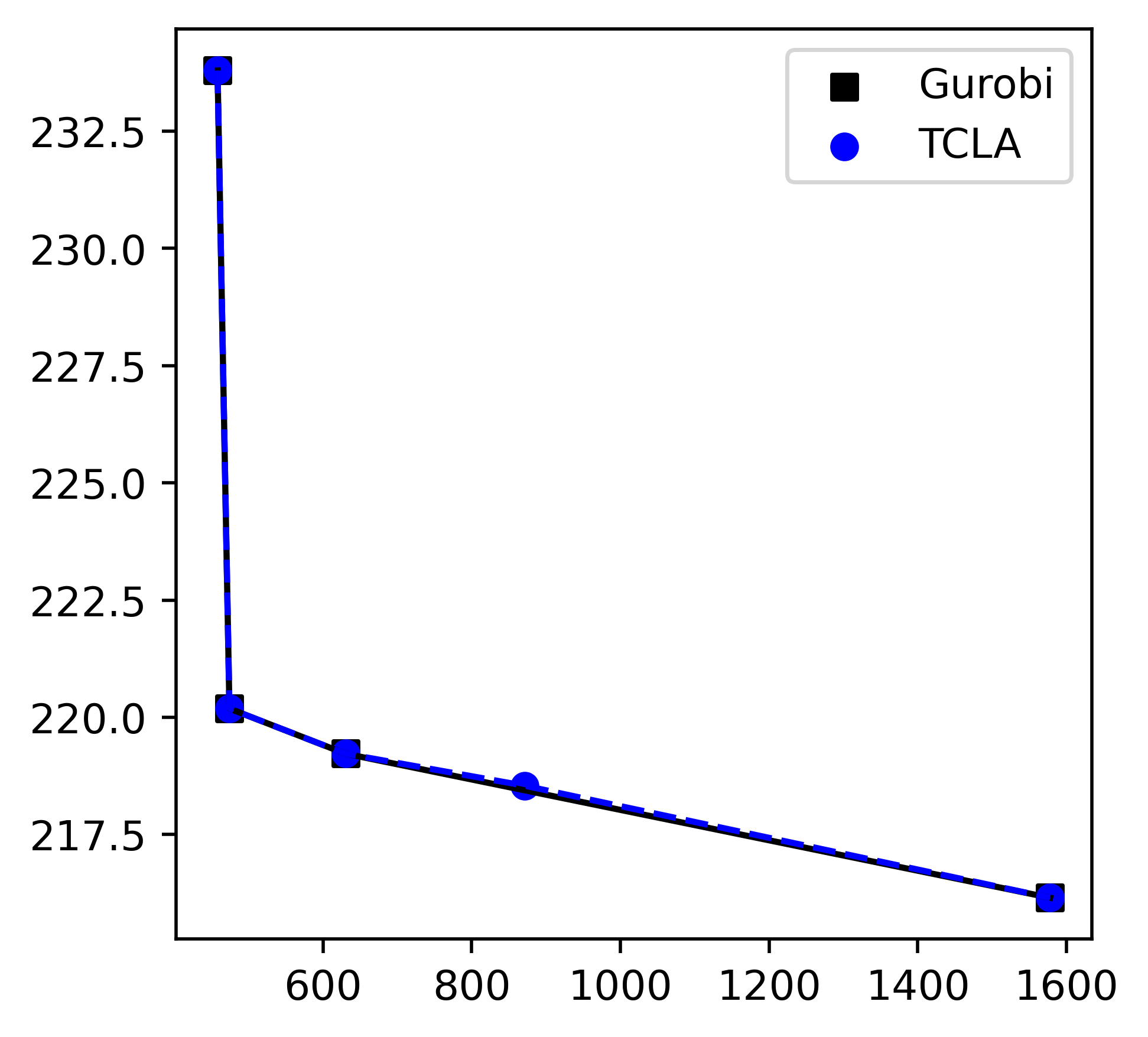}
        \caption*{instance (200,20,5)}
    \end{subfigure}
    \begin{subfigure}{0.24\textwidth}
        \centering
        \includegraphics[width=\linewidth]{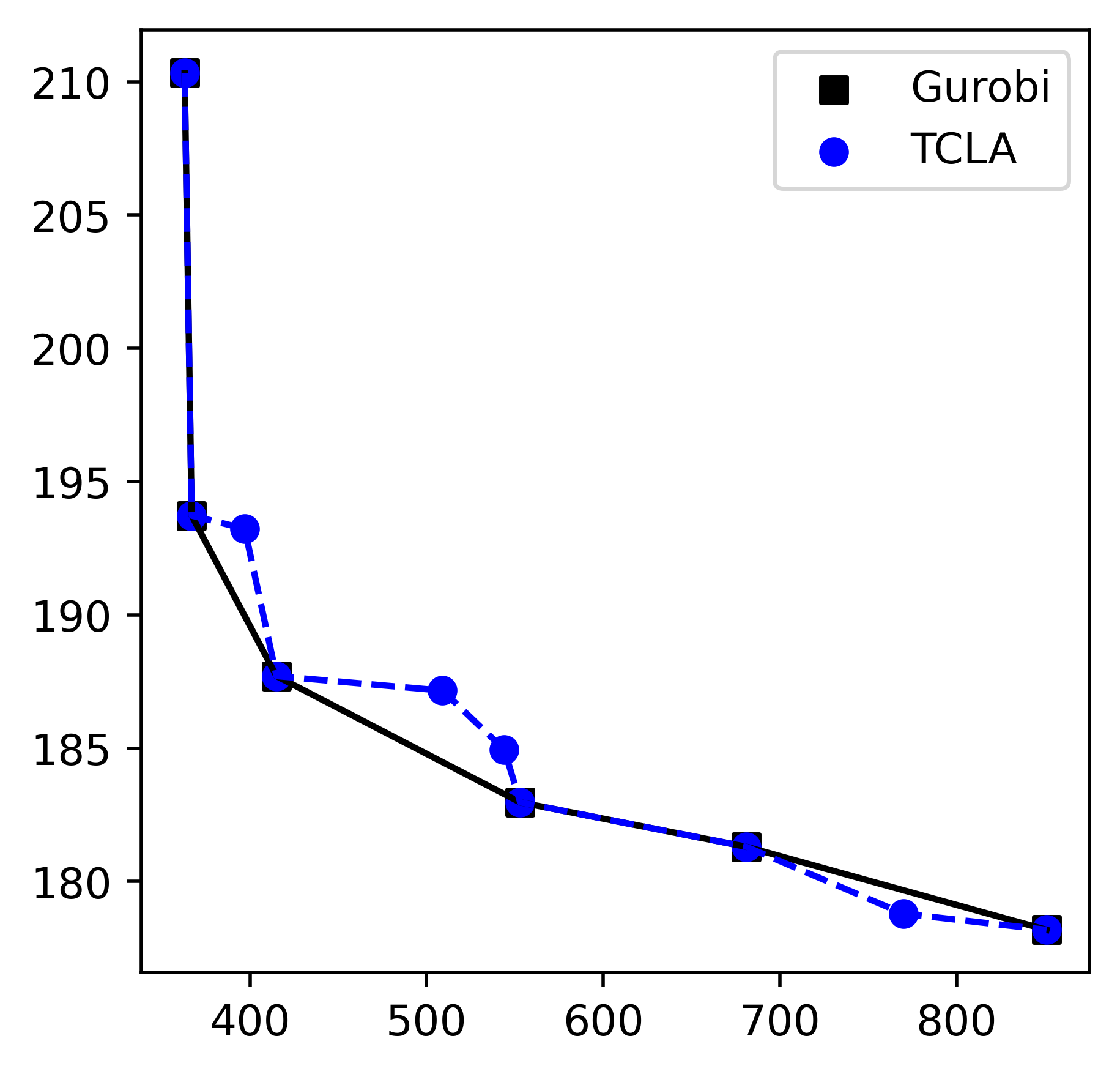}
        \caption*{instance (200,20,8)}
    \end{subfigure}
     \begin{subfigure}{0.24\textwidth}
        \centering
        \includegraphics[width=\linewidth]{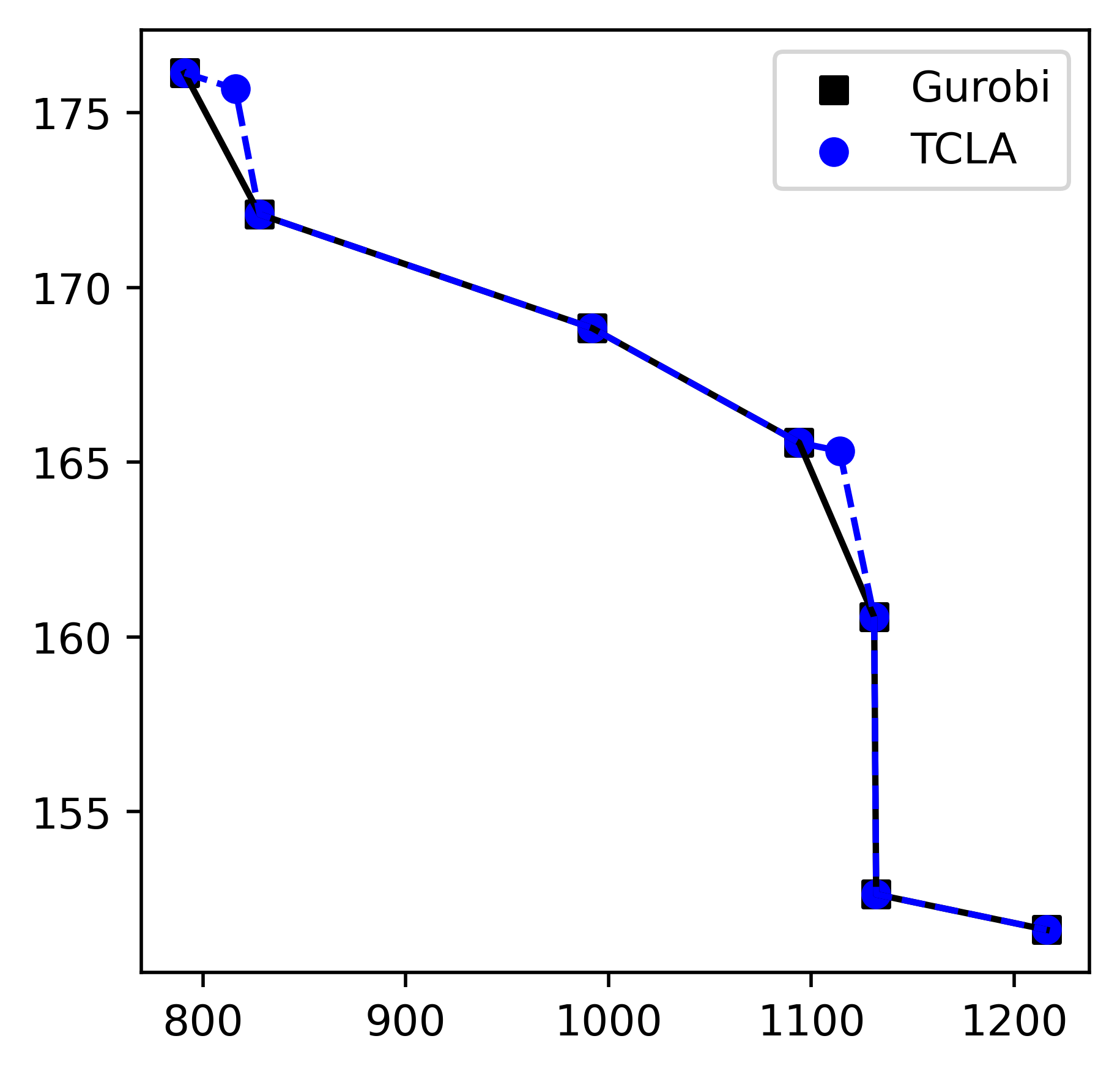}
        \caption*{instance (200,20,10)}
    \end{subfigure}
    \begin{subfigure}{0.24\textwidth}
        \centering
        \includegraphics[width=\linewidth]{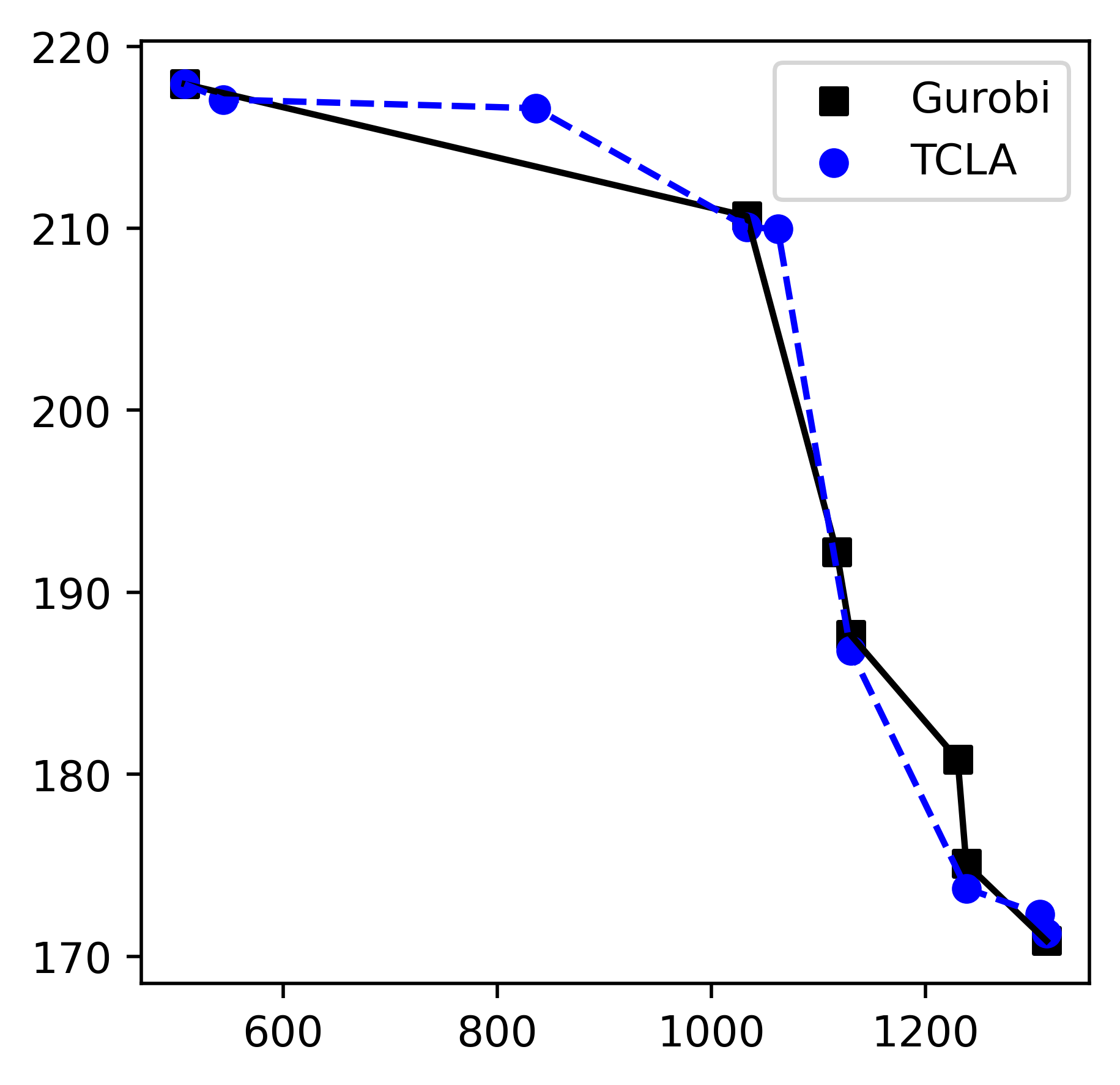}
        \caption*{instance (200,20,12)}
    \end{subfigure}

    \vspace{0.2cm} 

 \begin{subfigure}{0.24\textwidth}
        \centering
        \includegraphics[width=\linewidth]{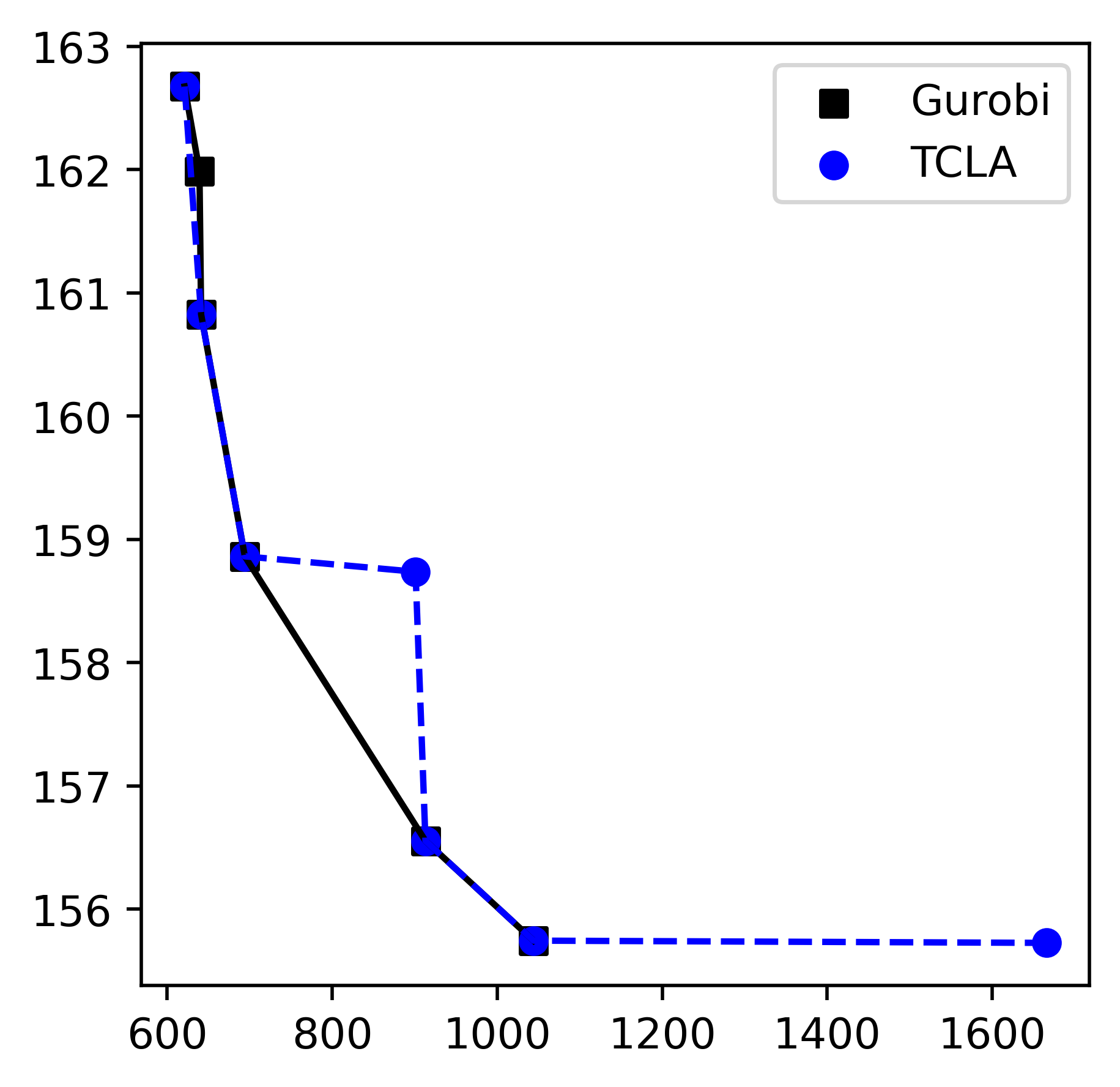}
        \caption*{instance (500,40,10)}
    \end{subfigure}
    \begin{subfigure}{0.24\textwidth}
        \centering
        \includegraphics[width=\linewidth]{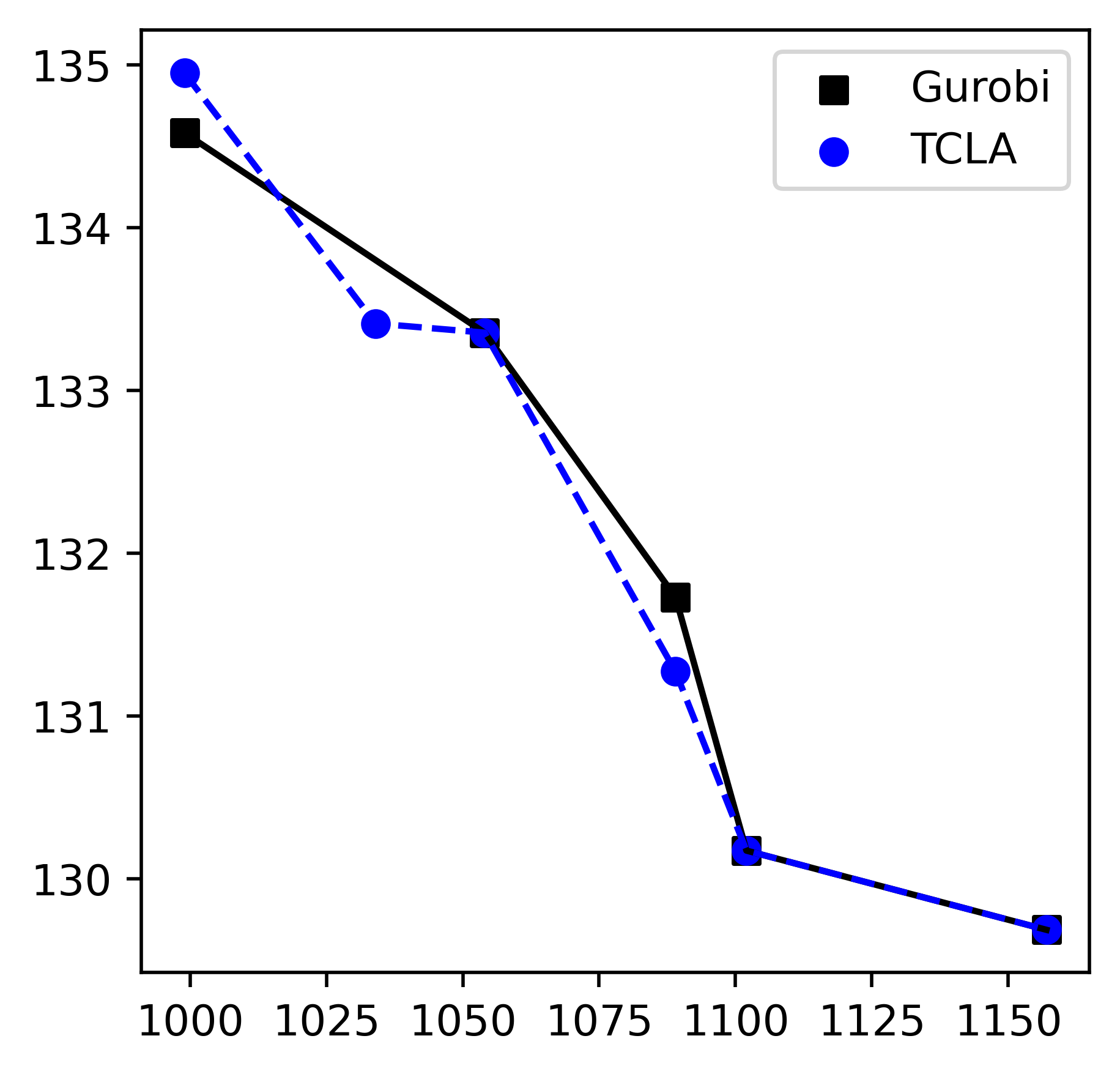}
        \caption*{instance (500,40,15)}
    \end{subfigure}
     \begin{subfigure}{0.24\textwidth}
        \centering
        \includegraphics[width=\linewidth]{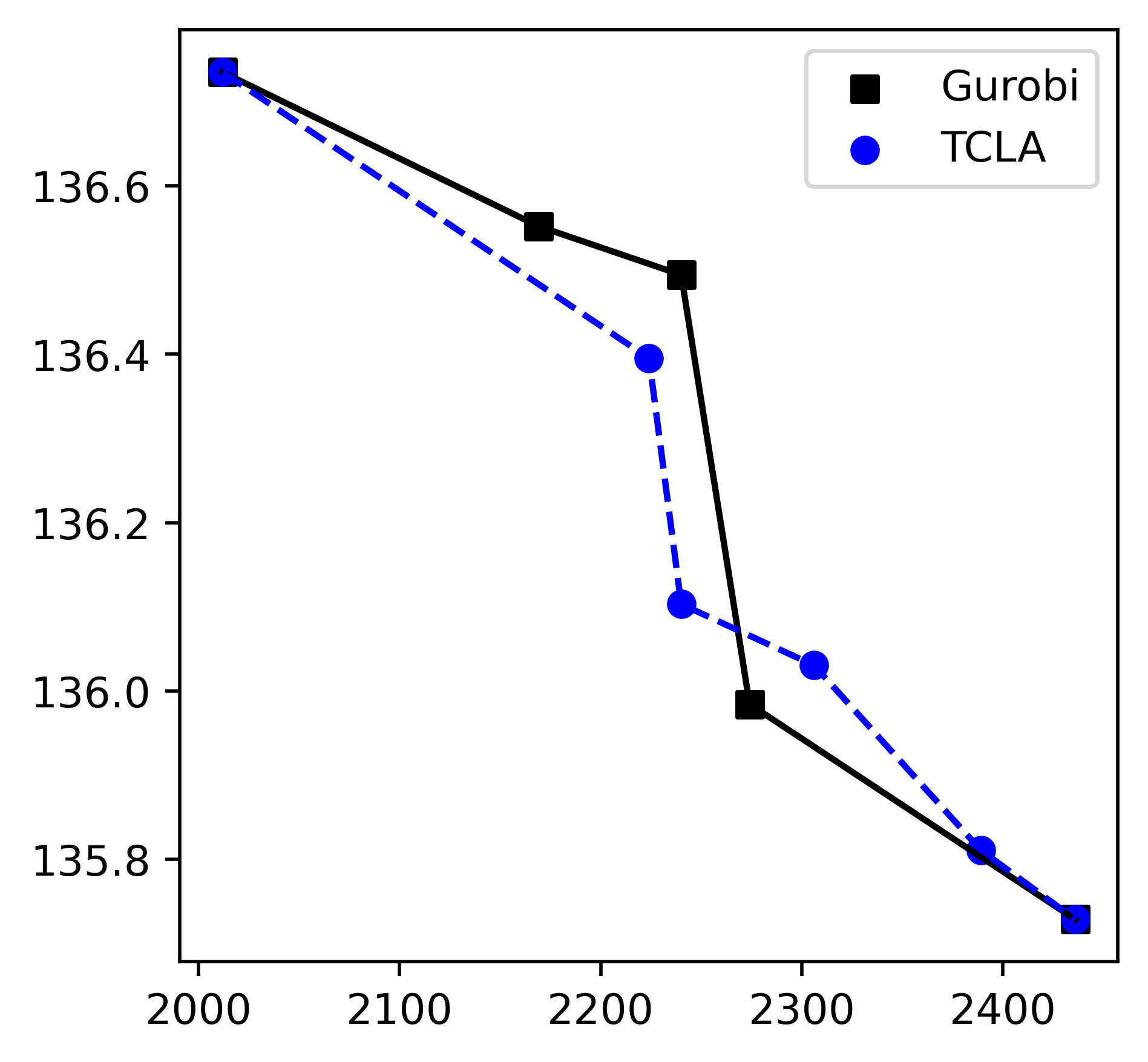}
        \caption*{instance (500,40,20)}
    \end{subfigure}
    \begin{subfigure}{0.24\textwidth}
        \centering
        \includegraphics[width=\linewidth]{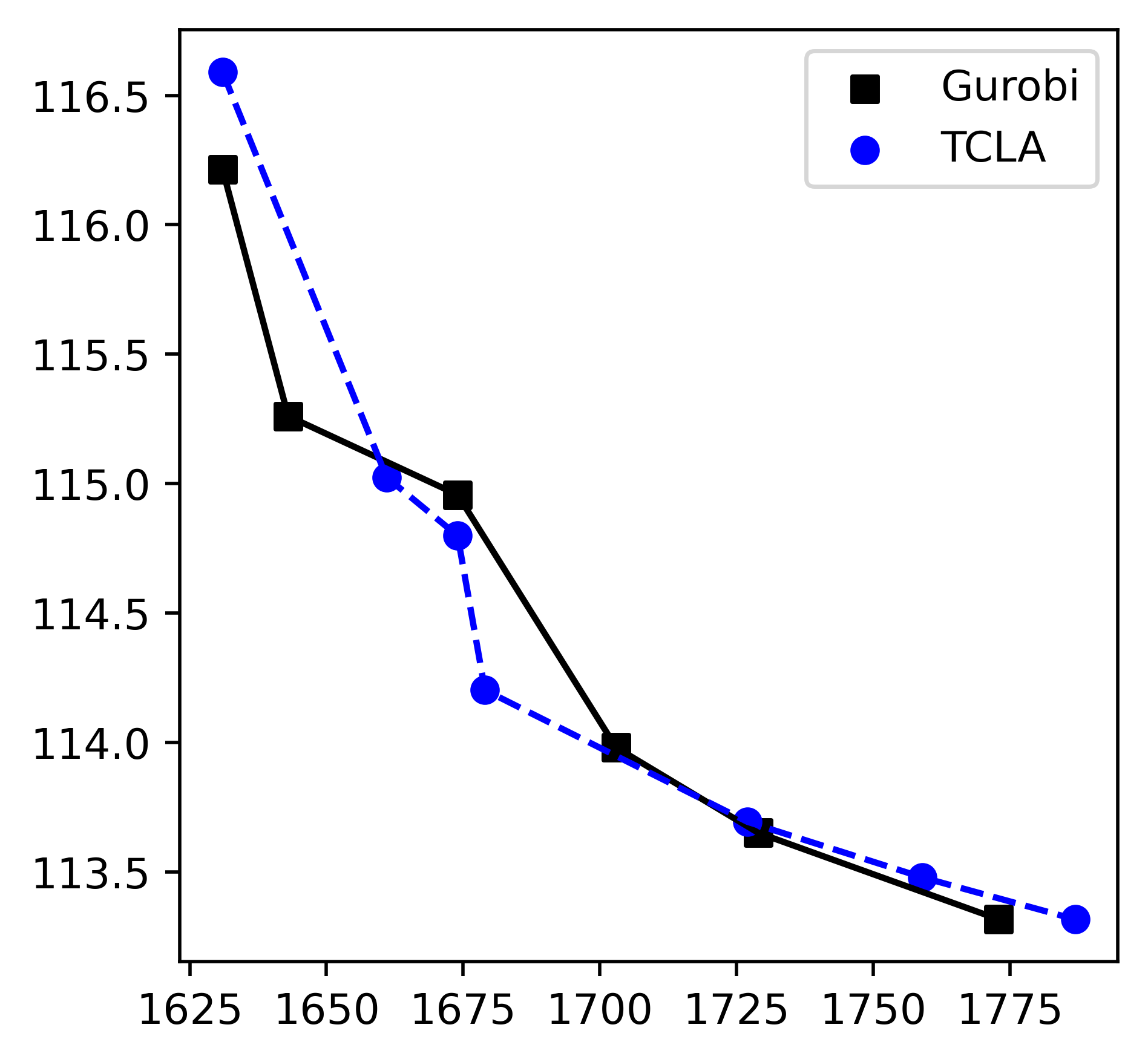}
        \caption*{instance (500,40,25)}
    \end{subfigure}

    \vspace{0.2cm} 
    
    \begin{subfigure}{0.49\textwidth}
        \centering
        \includegraphics[width=\linewidth]{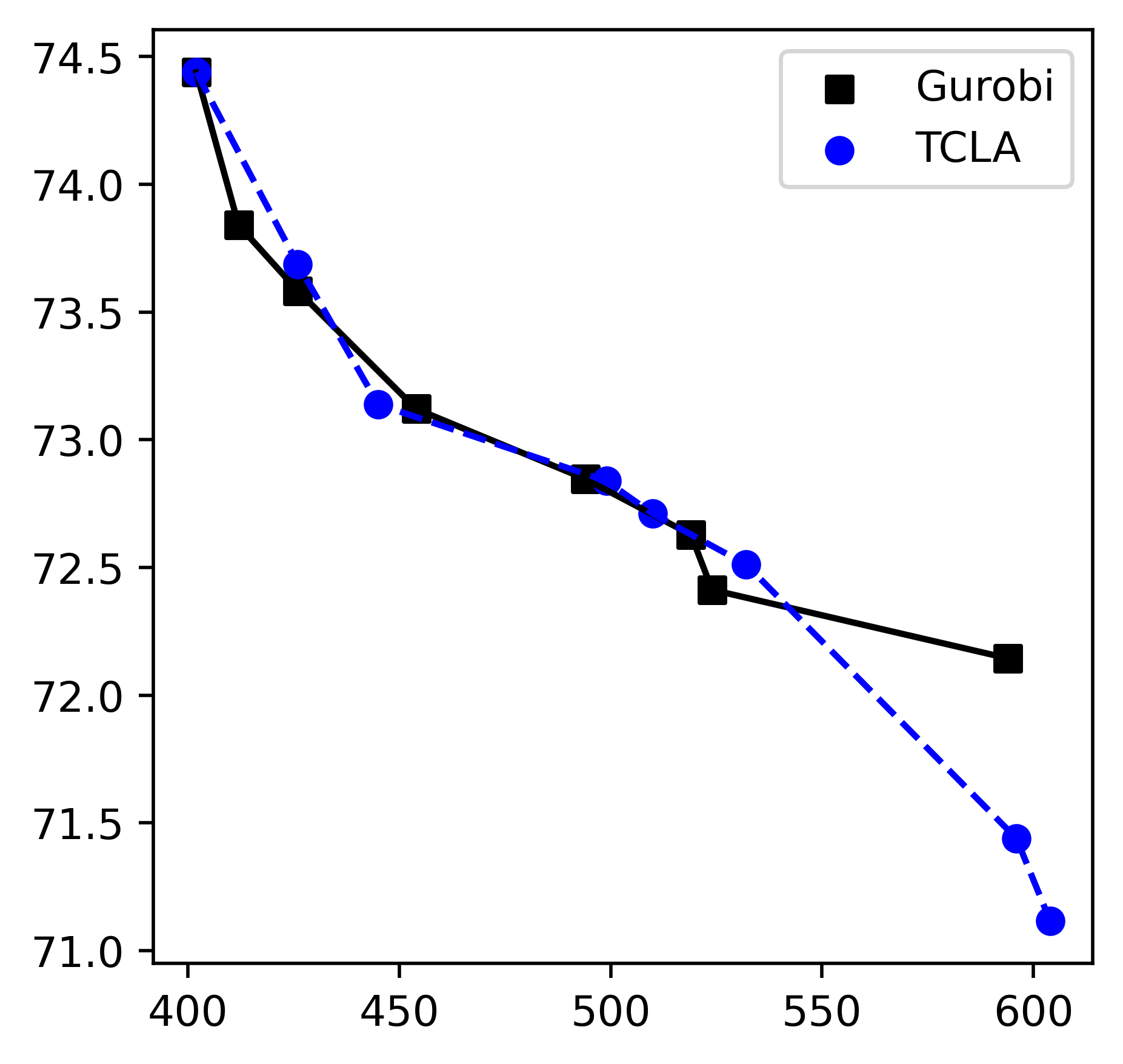}
        \caption*{instance (500,100,50)}
    \end{subfigure}
    \begin{subfigure}{0.49\textwidth}
        \centering
        \includegraphics[width=\linewidth]{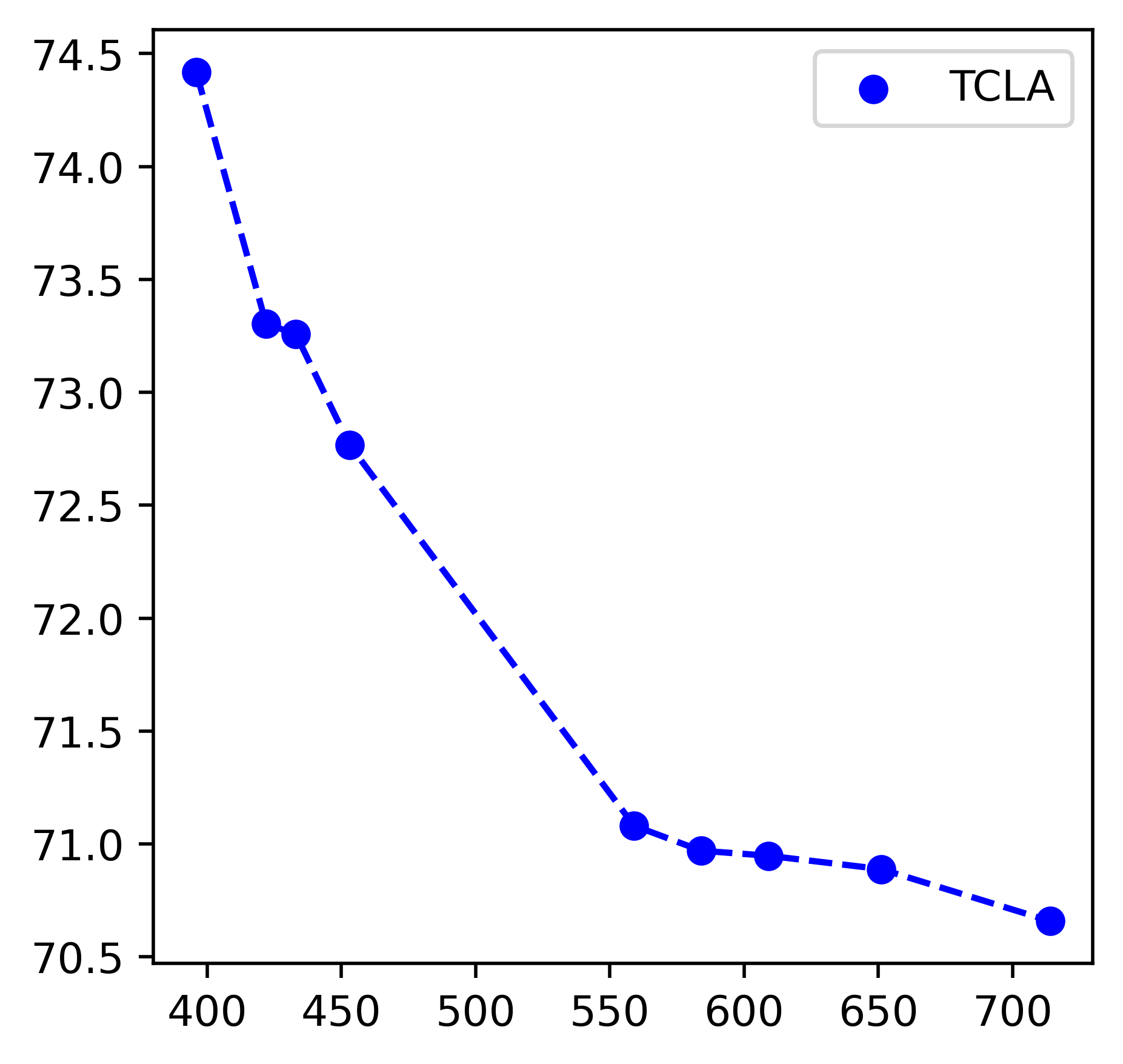}
        \caption*{instance (1000,100,50)}
    \end{subfigure}
    
    \caption{Obtained set by TCLA (blue-color diagram) and Gurobi (black-color diagram) in objective space for small instances. There is no outcome for Gurobi for the last instance.
    \label{ParetoTG}}
\end{figure}

The comparison results are presented in Table \ref{comparetable}, and the corresponding obtained non-dominated solutions in the objective space are depicted in Figure \ref{ParetoTG}. We ran the algorithms and asked to find at most 10 non-dominated solutions. The table represents the number of non-dominated solutions by each of the approaches (TCLA is denoted by \textit{T} and Gurobi is denoted by \textit{G}) as well as $scm(.,.)$ metric, $\alpha\beta(.,.)$ metric and the running time (in seconds) for all 18 different instances of the problem. The first 8 instances are generated randomly in a 150x100 rectangular shape, while for the larger instances with 200, 500, and 1000 weighted demand points, a bigger rectangle with a size of 1500x1000 is used. The weights are also assigned randomly in the interval $[10,100]$. The final two runs, pertaining to the instances (500,100,50) and (1000,100,50), represent exceedingly large cases of the TCLP. In these scenarios, there exists a staggering number of possible center combinations, approximately on the order of $5.39 \times 10^{23}$.
To conduct TCLA runs for these instances, we configured the population size to $N=1000$ and the number of generations to $T=15000$. The computational time for TCLA under these settings amounted to 1061 seconds for the (500,100,50) instance and 1937 seconds for the (1000,100,50) instance.
In contrast, Gurobi encounters substantial challenges when dealing with such formidable instances. For the former case, it necessitates a staggering 12841 seconds (over 3 hours and 30 minutes), while for the latter, we were compelled to terminate the program after 6 hours due to a lack of any discernible outcome. It's worth noting that TCLA exhibited a relatively modest memory usage of approximately 200 MB. In contrast, when using Gurobi, the memory consumption significantly surpassed this, exceeding 3800 MB. Moreover, the ensuing results have been derived from the comparison between TCLA and Gurobi.

\begin{table}[h]
\centering
\caption{Comparison results for TCLA (denoted by $T$) and Gurobi (denoted by $G$). Columns $\#T$ and $\#G$ show the number of obtained solutions by each method. The $scm$ and $\alpha \beta$ metrics are shown in the middle columns, and finally, the running time of the algorithms is reported in the last columns. For simplicity, $\alpha \beta(.,.)$ values are represented as percentage (*100).}
\scriptsize
\begin{tabular}{lcccccccc}
\hline
\textit{\textbf{(n,m,k)}} & \textit{\textbf{$\#T$}} & \textit{\textbf{$\#G$}} & \textit{\textbf{scm(T,G})} & \textit{\textbf{scm(G,T)}} & \multicolumn{1}{c}{\textit{\textbf{$\alpha\beta(T,G)$}}} & \multicolumn{1}{c}{\textbf{$\alpha\beta(G,T)$}} & \textit{\textbf{t(T)}} & \textit{\textbf{t(G)}} \\ \hline
(40,20,5)        & 4            & 4            & 0                 & 0                 & (0,0)                              & (0,0)                     & 3.81          & 14.48         \\
(40,20,8)        & 8            & 7            & 0                 & 0                 & (0,0)                              & (0,0)                     & 13.53         & 15.7          \\
(40,20,10)       & 10           & 8            & 0.25              & 0                 & (0,1.31)                           & (0,0)                     & 23.47         & 22.68         \\
(40,20,12)       & 7            & 6            & 0                 & 0                 & (0,0)                              & (0,0)                     & 12.47         & 11.38         \\
(100,40,10)      & 10           & 8            & 0                 & 0               & (0,0)                              & (0,0)               & 135           & 151           \\
(100,40,15)      & 8           & 9            & 0                 & 0.11               & (0,0)                              & (2.23,0.21)               & 207           & 105           \\
(100,40,20)      & 5            & 5            & 0.2              & 0.2                 & (3.67,0.85)                           & (0,0.35)                     & 237           & 114           \\
(100,40,25)      & 8            & 9            & 0.22              & 0.125             & (0,0.12)                           & (0,0.12)                  & 216           & 88.9          \\
(200,20,5)       & 5            & 4            & 0                 & 0                 & (0,0)                              & (0,0)                     & 17.98         & 89.06         \\
(200,20,8)       & 10           & 6            & 0                 & 0                 & (0,0)                              & (0,0)                     & 64.68         & 101.9         \\
(200,20,10)      & 9            & 7            & 0                 & 0                 & (0,0)                              & (0,0)                     & 96.29         & 109.5         \\
(200,20,12)      & 9            & 7            & 0.43              & 0.11              & (0,0.77)                           & (0,0.24)                  & 95.76         & 123.5         \\
(500,40,10)      & 7            & 6            & 0                 & 0                 & (0,0)                              & (0,0)                     & 277.4         & 1494          \\
(500,40,15)      & 6            & 5            & 0.2               & 0.17              & (0,0.35)                           & (0,0.27)                  & 576           & 1086          \\
(500,40,20)      & 6            & 5            & 0.2               & 0.17              & (0.71,0.29)                        & (1.39,0.03)               & 727           & 1884          \\
(500,40,25)      & 7            & 6            & 0.17              & 0.26              & (0,0.14)                           & (0.78,0.32)               & 637           & 1452          \\  \hline
(500,100,50)      & 8            & 8            & 0              & 0.25              & (0,0)                           & (1.50,0.14)               & 1061           & 12841          \\ 
(1000,100,50)      & 9            & --            & --              & --              & --                           & --               & 1737           & --          \\ \hline
\end{tabular}
\label{comparetable}
\end{table}

\begin{itemize}
    \item A trade-off exists between the running times of both algorithms and the quality of the non-dominated solutions they produce. As previously mentioned, the running time of TCLA is directly influenced by the size of potential center locations, denoted as $m$, and the number of selected centers, denoted as $k$. To simplify its application, we set the population size to $N = 2cm$ and the number of generations to $T = cN$, where $c$ is the minimum of either $k$ or $m-k$.
    On the other hand, in the case of Gurobi, since the model is an integer linear program, its solution quality is influenced by the parameter $MIPGap$. Smaller values of $MIPGap$ result in higher solution quality but longer running times. 
    In general, increasing the number of demand points has a noticeable impact on the running time of both Gurobi and TCLA. However, the impact is less pronounced in the case of TCLA. For larger instances, Gurobi may take more than 3 hours to complete.
    \item Another important consideration is that the running time of Gurobi is closely tied to the spatial distribution of points and the shape of the search space. In addition to instance size, the specific locations of the points play a crucial role in Gurobi's performance. However, this factor has a less significant impact on the running time of TCLA.
    \item The average \textit{scm} metric for $scm(TCLA,Gurobi)$ and $scm(Gurobi,TCLA)$ are 0.098 and 0.081, respectively. So, TCLA finds solutions slightly close to real Pareto-optimal solutions compared to Gurobi's solutions. However, because of the small size of the obtained non-dominated solutions, it appears that the $\alpha\beta$ metric provides a better assessment of solution quality.
    \item According to $\alpha\beta$ metric, aside from two instances, namely $(100,40,15)$ and $(100,40,20)$, where the difference between the obtained solutions is approximately 2.23\% and 3.67\%, respectively, the efficiency of both approaches is similar.
    \item While the Pareto-optimal queue of the TCLP is generally non-convex, the $\epsilon$-constraint approach demonstrates an ability to discover a diverse array of solutions. 
    \item The reported time for the Gurobi is the total time to run the Gurobi for different $\epsilon$ values in the $\epsilon$-constraint approach. So, if a decision maker is interested in finding just one single Pareto solution with a preferable level of work balance or average travel distance, the Gurobi will run faster than TCLA on small and medium-sized instances. 

\end{itemize}

Generally, the comparison results indicate that both TCLA and the Gurobi-based $\epsilon$-constraint approach show significant promise in effectively tackling the test center location problem. The choice between these approaches may depend on various factors such as problem size and computational resources (time and space), with each approach demonstrating its advantages.
TCLA excels in providing swift and reasonably high-quality solution sets, making it particularly suitable for scenarios where quick decision-making is essential. On the other hand, the $\epsilon$-constraint approach with the Gurobi solver offers a quicker solution for  small and medium-sized instances of the TCLP, especially when the objective is to identify a single optimal solution. This advantage stems from Gurobi's proficiency in handling integer linear programming models, while TCLA proves its adaptability in situations where linearity is not a critical constraint. As a result, TCLA holds the potential for broader applicability and extension to various problem variations, especially those demanding nonlinear modeling, such as scenarios where distance calculations, such as Euclidean distance, should be integrated directly into the model.
Finally, it becomes evident that Gurobi struggles to handle large-size instances of the TCLP within a reasonable time frame, whereas TCLA demonstrates competent performance and successfully identifies acceptable non-dominated solutions. One way to tackle this issue is by improving the formulation in Eq. (\ref{Model}). For example, Marín \cite{marin2011discrete} added some valid inequalities that help efficient branching and pruning in the branch-and-bound algorithms. Unfortunately, such improved formulation works only for unweighted demand points and it is not straightforward to extend this approach to the TCLP presented in this paper.


\section{Conclusion}
\label{conclusion}
In this paper, we have tackled a critical concern related to the establishment of diagnostic test centers for infectious diseases, drawing inspiration from the testing capacity limitations repeatedly exposed during the COVID-19 pandemic. Our primary objectives were to reduce workload disparities among centers while concurrently minimizing the average travel distance for individuals seeking testing. This posed a multifaceted challenge with significant real-world implications.
To address this complex problem, we introduced an integer linear programming model. Additionally, we proposed two distinct approaches for its solution. The first is a local search algorithm, named TCLA, which leverages Voronoi diagrams to efficiently uncover a set of non-dominated solutions in a single execution. The second approach employs the $\epsilon$-constraint method, solved using the Gurobi solver. We conducted comprehensive testing across a range of problem instances, rigorously assessing the performance of these approaches in terms of computational time and the quality of the resultant non-dominated solutions. Here, quality is gauged by the proximity of the obtained solutions to the Pareto-optimal solutions.
In light of the trade-off between computational time and solution quality, our comparative analysis demonstrates that TCLA emerges as an efficient algorithm for identifying Pareto-optimal solutions within a reasonable timeframe. This efficiency is particularly evident in the context of larger problem instances, where TCLA  outperforms Gurobi. This suggests its practical utility in real-world scenarios where time constraints are critical.

The models and approaches presented in this paper hold practical significance across a spectrum of real-world applications extending beyond healthcare systems. These principles can be applied to a wider range of facility location challenges where achieving workload equilibrium among centers is of utmost importance. Furthermore, by integrating elements such as uncertainties related to demand fluctuations or variations in travel times, as well as leveraging geographic information system data and spatial analysis, it is possible to create more realistic models that better align with real-world scenarios. 
Furthermore, in certain scenarios, the feasible facility center locations can be continuous, allowing for the possibility of opening centers in various positions throughout the city. For instance, during the COVID-19 pandemic, small kiosks offered antigen tests, illustrating this flexibility. In such cases, the model presented in this study may not be applicable, necessitating the development of a new formulation.

\bibliography{bibfile}
\bibliographystyle{elsarticle-harv}

\end{document}